\documentclass[11pt]{article}
\usepackage{latexsym, amsmath,amsbsy}
\newcommand{\nref}[1]{(\ref{#1})}
\def\Cal{\cal}

\newtheorem{theorem}{Theorem}
\newtheorem{remark}{Remark}[section]
\def\text#1{\hbox{#1}}

\def\endproof{\mbox{\ $\Box$}}
\def\Cal{\cal}

\newcommand{\CN}{{\Cal{N}}}
\newcommand{\CM}{{\Cal{M}}}
\newcommand{\Cl}{{\cal{H}}}

\def\a{\alpha}
\def\t{\theta}

\def\e{\varepsilon}
\def\b{\beta}
\def\g{\gamma}
\def\de{\delta}
\def\la{\lambda}

\def\d_1{\gamma_1}
\def\l{\left}
\def\r{\right}

\def\Var{\rm {Var}}
\def\tr{{\rm {Tr}}}

\def\CU{{\cal{U}}}
\def\CF{{\cal{F}}}

\def\CL{{\cal{L}}}

\def\Var{{\rm Var}}

\newtheorem{proposition}{Proposition}
\newtheorem{lemma}{Lemma}[section]

\newtheorem{corollary}{Corollary}

\newtheorem{definition}{Definition}[section]

\newcommand{\eq}        {{ \,\stackrel{\Delta}{=}\, }}
\newcommand{\iid}        {{ \,\stackrel{i i d}{\sim}\, }}

\newcommand{\N}        {{{\rm I\! N}}}
\newcommand{\R}        {{{\rm I\! R}}}

\newcommand{\Z}        {{{\rm Z\!\! Z}}}

\def\1{\mbox{1\hspace{-.30em}I}}
\numberwithin{equation}{section}

\begin{document}

\title{Minimax Goodness-of-Fit Testing in Multivariate Nonparametric Regression}

\author{Yuri I. Ingster \footnote{Research was partially supported  by RFBI
Grant 08-01-00692-a and by Grant NSh--638.2008.1}
\\ {\small Department of Mathematics II,} \\
{\small St. Petersburg State Electrotechnical University,}\\
{\small  Russia}\\ \\ and\\ \\Theofanis Sapatinas\\ {\small
Department of Mathematics and
Statistics,}\\
{\small  University of Cyprus,}\\
{\small  Cyprus}}

\date{}

\maketitle

\begin{abstract}
We consider an unknown response function $f$ defined on
$\Delta=[0,1]^d$, $1\le d\le\infty$, taken at $n$ random uniform
design points and observed with Gaussian noise of known variance. Given a positive
sequence $r_n\to 0$ as $n\to\infty$ and a known function $f_0 \in
L_2(\Delta)$, we propose, under general conditions, a unified
framework for the goodness-of-fit testing problem for testing the
null hypothesis $H_0: f=f_0$ against the alternative $H_1:
f\in\CF,\, \|f-f_0\|\ge r_n$, where $\CF$ is an ellipsoid in the
Hilbert space $ L_2(\Delta)$ with respect to the tensor product
Fourier basis and $\|\cdot\|$ is the norm in $ L_2(\Delta)$. We
obtain both rate and sharp asymptotics for the error probabilities
in the minimax setup. The derived tests are inherently non-adaptive.

Several illustrative examples are presented. In particular, we
consider functions belonging to ellipsoids arising from the
well-known multidimensional Sobolev and tensor product Sobolev norms
as well as from the less-known Sloan-Wo$\rm\acute{z}$niakowski norm
and a norm constructed from multivariable analytic functions on the
complex strip.

Some extensions of the suggested minimax goodness-of-fit testing
methodology, covering  the cases of general design schemes with a known product probability density function, unknown variance, other basis functions and adaptivity of the suggested tests, are also briefly discussed.

\end{abstract}


\medskip
\medskip


\newpage

\section{Introduction}
We consider the multivariate nonparametric regression model with a
random uniform design. More precisely, we observe
\begin{equation}
\label{eq:nprF} x_i=f(t_i)+\xi_i,\quad i=1,\ldots,n,
\end{equation}
where $t_i$ are random design points, $t_i\in \Delta = [0,1]^d$, $1
\leq d \leq \infty$. In particular, we assume that $t_i=\{t^k_i\}$
are (for $k=1,\ldots,d$ and $i=1,\ldots,n$) independent and identically distributed ({\em
iid}) random variables with a uniform distribution, i.e., $t^k_i
\iid \CU(0,1)$. Moreover, we assume that, conditionally on
$T_n=\{t_1,\ldots,t_n\}$, $\xi_i$ are {\em iid} Gaussian random variables with mean zero and variance $\tau^2$, i.e.,
$\xi_i\iid \CN(0,\tau^2)$, where
$\tau^2$ is assumed to be {\em known} with $0 < \tau^2 <
\infty$.

Given a positive sequence $r_n\to 0$ as $n\to\infty$ and a {\em known}
function $f_0 \in  L_2(\Delta)$, where $L_2(\Delta)$ is the set of
squared-integrable functions on $\Delta$, we propose, under general
conditions, a unified framework for the goodness-of-fit testing
problem for testing the null hypothesis
\begin{equation}
\label{eq:nullF} H_0: f=f_0
\end{equation}
against the alternative
\begin{equation}
\label{eq:altF} H_1: f\in\CF,\,\|f-f_0\|\ge r_n,
\end{equation}
where $\CF$ is an ellipsoid in the Hilbert space $ L_2(\Delta)$ with
respect to the tensor product Fourier basis and $\|\cdot\|$ is the
norm in $ L_2(\Delta)$. (The set $\CF$ corresponds to a ``regularity
constraint'' on the response function $f$.)

We are interested in both rate and sharp asymptotics for the error
probabilities in the minimax setup, i.e., we try to find the maximal
rate of convergence of $r_n\to 0$ as $n\to\infty$ which provide
nontrivial minimax testing, when certain constraints are imposed on
the regularity of the response function $f$.

Although there is a plethora of research work in the literature on
the estimation problem for response functions $f\in \CF$ in (both
univariate and multivariate) nonparametric regression (under various
design schemes), much less attention has been paid to the hypotheses
testing problem in this model, especially in the multivariate case.
This work is devoted to the goodness-of-fit testing problem
(\ref{eq:nullF})--(\ref{eq:altF}) in the nonparametric regression
model (\ref{eq:nprF}).

Nonparametric goodness-of-fit testing was studied intensively during
the last twenty years or so; however,  main results were
obtained for the detection of the response function $f \in
L_2(\Delta)$, with $d=1$, in the 1-variable Gaussian white noise
model, i.e.,
\begin{equation}\label{GM}
dX(t)=f(t)dt+\e dW(t),\quad t\in [0,1],
\end{equation}
where $W(t)$ is the standard Wiener process, with the noise level
$\e\to 0$. In particular, rate and sharp asymptotics for the error
probabilities in the minimax setup were obtained for various
classes $\CF$ of nonparametric alternatives. Moreover, under
periodicity, the sharp asymptotics are of Gaussian type and are
determined by a specific extremal problem (see, e.g., \cite{E.90a},
\cite{E.03}, \cite{I.93}, \cite{IS.02}).

These results have been extended in part to the density, spectral
density, nonparametric regression and  Poisson models for the
1-variable case (see, e.g., \cite{E.03}, \cite{I.93}, \cite{IK.07},
\cite{IS.02}). Note that, under some regularity constraints, one can
formally deduce some results for the 1-variable density and
nonparametric regression models from results on the asymptotic
equivalence (in Le Cam sense) of these models to the 1-variable
Gaussian white noise model (see, e.g., \cite{BL}, \cite{N}).

For the $d$-variable Gaussian white noise model, we have typically
similar separation rates with the smoothness parameter $\sigma$
(associated with the ``regularity constraint'' on the response
function $f$) replaced by $\tilde\sigma=\sigma/d$ as well as sharp
asymptotics of a similar type (see \cite{IS.05}). This leads to the
``curse of dimensionality" phenomenon when $d$ is large (see
\cite{IS.06}). It was recently shown that one can actually lift the
curse of dimensionality by using different type of regularity
constraints, which are determined by the so-called
``Sloan-Wo$\rm\acute{z}$niakowski'' norm (see \cite{IS.06}).
Although, analogously to the 1-variable case, one can formally
deduce, under some stronger regularity constraints, some results for
the multivariate nonparametric regression models from results on the
asymptotic equivalence (in Le Cam sense) of these models to the
$d$-variable Gaussian white noise model (see, e.g., \cite{C.2006},
\cite{R}), one cannot apply these results to the tensor product
Sobolev or Sloan-Wo$\rm\acute{z}$niakowski type spaces, because
there are no asymptotic equivalence results as yet for these spaces.

Rate asymptotics in $d$-variable parametric regression models were
studied in, e.g., \cite{GL}, \cite{HS}, for testing a parametric
model against Lipschitz and H\"older classes $\CF$ of alternatives,
respectively. On the other hand, rate asymptotics in the
multivariate regression model, under equispaced design points, were
studied in \cite{ADFS} for the goodness-of-fit testing problem
(\ref{eq:nullF})--(\ref{eq:altF}), under Besov balls $\CF$ of
alternatives.

The purpose of this paper is to extend some results on the
goodness-of-fit testing  of \cite{E.90a}, \cite{I.93},
\cite{IS.02}-\cite{IS.07} for the $d$-variable Gaussian white noise
model to the goodness-of-fit testing problem
(\ref{eq:nullF})--(\ref{eq:altF}) for the multivariate nonparametric
regression model (\ref{eq:nprF}), in a unified framework.

In our study, we use analytic results on an extermal problem
for ellipsoids that were presented in \cite{I.93},
\cite{IS.02}-\cite{IS.07} for the $d$-variable Gaussian white noise
model. These lead to the asymptotic efficiency of testing  for the
multivariate nonparametric regression model (\ref{eq:nprF}), similar
to the ones that have earlier been obtained, in specific settings,
for the $d$-variable Gaussian white noise model, under the standard
calibration $\varepsilon=\tau/\sqrt{n}$. However, the machinery of
reduction of the hypothesis testing problems to the extermal problem is
different and, essentially, more difficult, especially for the study of the lower
bounds. The proposed tests are of different structure as well: they are
based on U-statistics of increasing dimension. Certainly, this
reduction requires some assumptions on the basis functions and on
the sample size (compare with \cite{Efr} for estimation problem). It
is a typical situation for extending results from the Gaussian white noise model to other statistical models (e.g., density, spectral
density, intensity of a Poisson process and so
on).

Several illustrative examples are presented. In particular, we
consider functions belonging to the balls under the well-known
multidimensional Sobolev and tensor product Sobolev norms as well as
from the less-known Sloan-Wo$\rm\acute{z}$niakowski norm and a norm
constructed from multivariable analytic functions on the complex
strip. Some extensions of the suggested minimax goodness-of-fit
testing methodology, covering  the cases of general design schemes with a known product probability density function, unknown variance, other basis functions and adaptivity of the suggested tests, are also briefly discussed.

\section{Preliminaries and assumptions}

\subsection{Minimax goodness-of-fit testing}
Consider the multivariate nonparametric regression model
(\ref{eq:nprF}). Given a known function $f_0 \in L_2(\Delta)$, we
test the null hypothesis (\ref{eq:nullF}), i.e., we test $H_0:
f=f_0$.  Given a positive sequence $r_n\to 0$ as $n \rightarrow
\infty$, let
$$
\CF(r_n)=\{f\in\CF: \|f-f_0\|\geq r_n\},
$$
where $\CF$ is an ellipsoid in the Hilbert space $ L_2(\Delta)$ with
respect to the tensor product Fourier basis and $\|\cdot\|$ is the
norm in $ L_2(\Delta)$. Consider now the alternative hypothesis
(\ref{eq:altF}), i.e., consider $ H_1: f\in\CF(r_n).$ (In what
follows, without loss of generality, we restrict ourselves to the
cases $f_0 = 0$ and $\tau=1$.)

Set $X_n=\{x_1,\ldots,x_n\}$ and recall that
$T_n=\{t_1,\ldots,t_n\}$. Let $P_{n,f}$ be the probability measure
that corresponds to $Z_n=(X_n,T_n)$ and denote by $E_{n,f}$ the
expectation over this probability measure. Let $\psi$ be a
(randomized) test, i.e., a measurable function of the observation
$Z_n$ taking values in $[0,1]$: the null hypothesis is rejected with
probability $\psi(Z_n)$ and is accepted with probability
$1-\psi(Z_n)$. Let
$$\a(\psi) = E_{n,0}\psi$$ be its type I error probability, and let
$$
\b(\CF,r_n,\psi) = \sup_{f\in \CF(r_n)}E_{n,f}(1-\psi)
$$
be its maximal type II error probability. We consider two criteria
of asymptotic optimality:

[1] The first one corresponds to the classical Neyman-Pearson
criterion. For $\a\in (0,1)$ we set
$$
\b(\CF,r_n,\a) = \inf_{\psi:\, a(\psi)\le\a}\b(\CF,r_n,\psi).
$$
We call a sequence of tests $\psi_{n,\a}$ {\it asymptotically
minimax} if
$$
\a(\psi_{n,\a})\le\a+o(1),\quad
\b(\CF,r_n,\psi_{n,\a}) = \b(\CF,r_n,\a)+o(1),
$$
where $o(1)$ is a sequence tending to zero; here, and in what
follows, unless otherwise stated, all limits are taken as $n
\rightarrow \infty$.

[2] The second one corresponds to the total error probabilities. Let
$\g(\CF,r_n,\psi)$ be the sum of the type I and the maximal type II
error probabilities, and let $\g(\CF,r_n)$ be the minimax total
error probability, i.e.,
$$
\g(\CF,r_n) = \inf_{\psi}\g(\CF,r_n,\psi),
$$
where the infimum is taken over all possible tests. We call a
sequence of tests $\psi_{n}$ {\it asymptotically minimax} if
$$
\g(\CF,r_n,\psi_n)=\g(\CF,r_n)+o(1).
$$
It is known that (see, e.g., Chapter 2 of \cite{IS.02}) that
$$
\b(\CF,r_n,\a)\in [0,1-\a], \quad  \g(\CF,r_n)=\inf_{\a\in
(0,1)}(\a+\b(\CF,r_n,\a)) \in [0,1].
$$

We consider the problems of rate and sharp asymptotics for the error
probabilities in the minimax setup. The rate optimality problem
corresponds to the study of the conditions for which $\g(\CF,r_n)\to
1$ and $\g(\CF,r_n)\to 0$ and, under the conditions of the last
relation, to the construction of asymptotically {\it minimax
consistent} sequences $\psi_n$, i.e, such that
$\g(\CF,r_n,\psi_n)\to 0$. Often, these conditions correspond to
some minimal decreasing rates for the sequence $r_n$. Namely, we say
that the positive sequence $r_n^*=r_n^*(\CF)$,  $r_n^* \rightarrow 0$, is a {\it separation rate}, if
$$
\g(\CF,r_n)\to 1\quad\text{as}\quad r_n/r_n^*\to 0,
$$
and
$$
\g(\CF,r_n)\to 0, \quad\text{and}\ \ \b(\CF,r_n,\a)\to 0 \ \
\text{for any}\ \a\in (0,1),\quad \text{as}\quad r_n/r_n^*\to
\infty.
$$
In other words, it means that, for large $n$, one can detect all
functions in $f \in \CF$ if the ratio $r_n/r_n^*$ is large, whereas
if this ratio is small then it is impossible to distinguish between
the null and the alternative hypothesis, with small minimax total
error probability. Hence, the rate optimality problem corresponds to
finding the separation rates $r_n^*$ and to constructing
asymptotically minimax consistent sequence of tests.

On the other hand, the sharp optimality problem corresponds to the
study of the asymptotics of the quantities $\b(\CF,r_n,\a),\
\g(\CF,r_n)$ (up to vanishing terms) and to the construction of
asymptotically minimax sequences $\psi_{n,\a}$, $\psi_{n}$,
respectively. Often, the sharp asymptotics are of Gaussian type,
i.e.,
\begin{equation}\label{G}
\b(\CF,r_n,\a)=\Phi(H^{(\a)}-u_n)+o(1),\quad
\g(\CF,r_n)=2\Phi(-u_n)+o(1),
\end{equation}
where $\Phi$ is the standard Gaussian distribution function,
$H^{(\a)}$ is its $(1-\a)$-quantile, i.e.,
$\Phi(H^{(\a)})=1-\alpha$, and the sequence $u_n=u_n(\CF,r_n)$
characterizes {\em distinguishability} in the problem. The
separation rates $r_n^*$ are usually determined by the relation
$u_n(\CF,r_n^*)\asymp 1$ (see, e.g.,  \cite{I.93}, \cite{IS.02}). Hence, the
sharp optimality problem corresponds to calculating the sequence
$u_n$ and to constructing asymptotically minimax sequence of tests.

\subsection{Assumptions}

Let $L_2(\Delta)=L_2$, $\CL$ be a denumerable set, $\{\phi_l\}_{
l\in\CL}$ be an orthonormal system in $L_2$, and $L_2^{\CL}\subset
L_2$ be the closed linear hull of the system $\{\phi_l\}_{l\in\CL}$.
For a function $f \in L_2^{\CL}$, let $\t=\{\t_l\}_{l\in\CL}$ be the
``generalized'' Fourier coefficients with respect to this system,
i.e., $\t_l =\langle f,\phi_l \rangle$, $l \in \CL$, where $\langle
\cdot,\cdot \rangle$ denotes the inner product in $L_2$.

Let a collection of coefficients $\{c_l\}_{l \in \CL}$, $c_l \geq0$,
be given. The set of functions $\CF \subset L_2^{\CL}$ under
consideration are the {\em ellipsoids} with respect to the
orthonormal system $\{\phi_l\}_{ l\in\CL}$ with coefficients
$\{c_l\}_{l \in \CL}$, ${l \in \CL}$, i.e.,
$$
\CF=\{f\ :\ f(t)=\sum_{l\in \CL}\t_l\phi_l(t),\ \sum_{l\in
\CL}c_l^2\t_l^2\le 1\}.
$$

Let
\begin{equation*}\label{NC}
\CN(C) = \{l\in \CL: c_l<C\}, \quad\ N(C) = \#\CN(C),
\end{equation*}
where $\#$ denotes the cardinality of a set.

\medskip

Consider the following set of assumptions:

\medskip

{\bf (A1)} The set $\CN(C)$ is finite, i.e.,
\begin{equation*}\label{NF}
 N(C)<\infty\quad \forall\ C>0.
\end{equation*}

{\bf (A2)} The orthonormal system $\{\phi_l\}_{ l\in\CL}$ satisfies
\begin{equation*}\label{Ph}
\sum_{l\in\CN(C)}\phi_l^2(t)
 = N(C) \quad \forall\ C>0,\
t\in\Delta.
\end{equation*}

\medskip
{{\bf (A3)} The functions $f \in \CF$ are uniformly bounded in
$L_p(\Delta)$-norm for some $p>4$, i.e.,
$$
\exists\ p>4:\quad \sup_{f\in\CF}\int_{\Delta}|f(t)|^p<\infty.
$$

\begin{remark}{\em
Note that assumption {\bf (A3)} follows from the following stronger
condition,
\begin{equation}\label{sup}
\sup_{f\in\CF}\|f\|_\infty<\infty,
\end{equation}
where $\|f\|_{\infty} = \sup_{t \in \Delta}|f(t)|$. }
\end{remark}

\section{Rate optimality}

In what follows, the relation $A_n \sim B_n$ means that $A_n/B_n$
tends to 1 while the relation $A_n \asymp B_n$ means that there
exists constants $0 <c_1 \leq c_2 <\infty$ and $n_0$ large
enough such that $c_1 \leq A_n/B_n \leq c_2$ for $n \geq n_0$. Let
also {\em $\1_{\{A\}}$} be the indicator function of a set $A$.

\medskip

For a sequence $C=C_n$, let $\CN = \CN(C_n)$, $N = N(C_n)$.

\medskip

Let us introduce an extra assumption.

\medskip
{\bf (B1)}\quad \quad \quad $N=o(n)$.

\begin{theorem}\label{T0} Let $r_n\to 0$.

(1)[Lower bounds] Assume {\bf (A1)}--{\bf (A2)}. Take $C_n\to\infty$
such that $\limsup(C_nr_n) < 1$ and {\bf (B1)} holds. Then
\begin{equation*}\label{G1}
\b(\CF,r_n,\a)\ge\Phi(H^{(\a)}-u_n)+o(1),\quad \g(\CF,r_n)\ge
2\Phi(-u_n)+o(1),
\end{equation*}
where
\begin{equation}\label{U1}
u_n^2=\frac{n^2r_n^4}{2N}.
\end{equation}

(2) [Upper bounds] Assume {\bf (A1)}--{\bf (A3)}. Take
$C_n\to\infty$ such that {\bf (B1)} holds. Consider the sequence of
tests $ \psi_{n}^{H}=\1_{\{U_n>H\}}$ based on the $U$-statistics
\begin{equation}
\label{eq:uF} U_n = \frac{1}{n}\sum_{1\le i<k\le n}K_n(z_i,z_k),
\end{equation}
where $z_i = (x_i,t_i),\ i=1,\ldots,n$ are the observations, with
the kernel
\begin{equation}\label{test1}
K_n(z^{'},z^{''}) = x^{'}x^{''}G_n(t^{'},t^{''}),\quad
G_n(t^{'},t^{''}) = \sqrt{\frac 2N} \sum_{l \in
\CN}\phi_l(t^{'})\phi_l(t^{''}).
\end{equation}
Set
\begin{equation}\label{test2}
h_n(f)=\frac{n}{\sqrt{2N}}\sum_{l\in \CN}\t_l^2.
\end{equation}
Then, uniformly over $H = H_n\in\R$,
$$
\a(\psi_{n}^{H}) \leq 1-\Phi(H)+o(1),
$$
and, for any $c\in (0,1)$, uniformly over $f\in\CF$ and $H = H_n$
such that $h_n(f)\ge cH_n$,
$$
\b(\CF,r_n,\psi_{n}^{H})\le \Phi(H-h_n(f))+o(1).
$$
\end{theorem}

\begin{remark} {\rm We now give some intuition about the suggested $U$-statistics
used in Theorem 1.
For testing the null hypothesis $H_0: f = 0$ in the Gaussian white
noise model, a natural test statistic is a centered and normalized
(under $H_0$) version of the quadratic functional $
\sum_{l\in\CL}\hat\t_l^2$, where $\hat\t_l =
\int_{\Delta}\phi_l(t)dX(t). $ The analog of $\hat\t_l$ in the
multivariate nonparametric regression model (\ref{eq:nprF}) is given
by $\hat\t_l = n^{-1}\sum_{i=1}^n\phi_l(t_i)x_i$ which leads to the
quadratic functional
$$
\sum_{l\in\CL}\hat\t_l^2=\frac{1}{n^2}\sum_{i,k=1}^nx_ix_k
\widetilde G_n(t_i,t_k),\quad \widetilde G_n(t^{'},t^{''}) =
\sum_{l\in\CL}\phi_l(t^{'})\phi_l(t^{''}).
$$
Suppressing now the terms with $i=k$, a centered and normalized
version of this quadratic functional corresponds to the
$U$-statistic defined in \nref{eq:uF} with the kernel defined in
\nref{test1}.}
\end{remark}

Let the sequence $C=C_n$ be determined by the ``balance equation''
\begin{equation}\label{balance}
C_n^4 N(C_n)\asymp n^2.
\end{equation}
Observe that, in this case, under {\bf (A1)}, $C_n\to\infty$ and,
hence, $\ N(C_n) \to \infty$.

\begin{remark} {\rm

Note that if $r_n $ satisfies $C_nr_n\asymp 1$, then (\ref{balance}) corresponds to $u_n\asymp 1$ in \nref{U1}. Corollaries \ref{C1a} and \ref{C2a} below show a motivation of
\nref{balance}. }
\end{remark}

\medskip Let us introduce an extra assumption.

\medskip

{\bf (B2)}\quad \quad For any $B>0$, $ N(C_n)\asymp N(BC_n)$.

\medskip

Note that we can obtain lower bounds for $h_n(f)$ from \nref{test2}.
Indeed, for $f\in\CF(r_n)$, we have
\begin{eqnarray}\nonumber
h_n(f)&=&\frac{n}{\sqrt{2N}}\bigg(\sum_{l\in\CL}\t_l^2-\sum_{c_l\ge
C_n}\t_l^2\bigg)\ge
\frac{n}{\sqrt{2N}}\bigg(r_n^2-C_n^{-2}\sum_{c_l\ge
C_n}c_l^2\t_l^2\bigg)\\ \label{test3}
&\ge&\frac{n}{\sqrt{2N}}(r_n^2-C_n^{-2})=\frac{nr_n^2}{\sqrt{2N}}\l(1-(r_nC_n)^{-2}\r).
\end{eqnarray}
Therefore, if $C_nr_n\ge B>1$, we have from Theorem \ref{T0} (2),
$$
\b(\CF,r_n,\psi_{n}^{H})\le \Phi\big(H-u_n(1-B^{-2})\big)+o(1),
$$
with $u_n$ determined by \nref{U1}. This leads to

\begin{corollary}\label{C1a} Let $r_n\to 0$. Assume {\bf (A1)}--{\bf (A3)} and {\bf (B1)}--{\bf (B2)}.
Then

[1] The separation rates are of the form
$$
r_n^* \asymp C_n^{-1},
$$
where the sequence $C =C_n$ is determined by \nref{balance}.

[2] Moreover, let $r_n/r_n^*\to\infty$. Then, there exists a
sequence $H=H_n\to\infty$ such that the sequence of tests
$\psi_{n}^{H}=\1_{\{U_n>H\}}$ is asymptotically minimax consistent,
i.e., $\g(\CF,r_n,\psi_{n}^{H})\to 0$.
\end{corollary}

We say that a function $g(t),\ t>0$, is a {\it slowly varying}
function if $g(Bt)/g(t)$ tends to 1 as $t\to\infty$, for any $B>0$.

\medskip
This leads to the following assumption.

\medskip

{\bf (B3)}\quad $N(C_n)$ is a  slowly varying function.

\begin{corollary}\label{C2a} Let $r_n\to 0$. Assume {\bf (A1)}--{\bf (A3)}   and {\bf (B1)}--{\bf
(B3)}. Then

[1] The sharp asymptotics \nref{G} hold, where $u_n$ is defined by
\nref{U1} with any $N(C_n)$ determined by \nref{balance}.

[2] Moreover, for any sequence $C_n$ satisfying \nref{balance},
there exists a sequence $B_n \to \infty$ such that, for the sequence
$C_{n,1}= B_nC_n$, the sequence of tests $\psi_{n}^{H^{(\a)}}$ is
asymptotically minimax under the Neyman-Pearson criterion, and the
sequence of tests $\psi_{n}^{u_n/2}$ is asymptotically minimax under
the total error probability criterion.
\end{corollary}
{\bf{Proof}}.
In order to get the upper bounds, note that under {\bf (B3)} one can
take a sequence $B_n\to\infty$ such that $N(B_nC_n)\sim N(C_n)$.
Applying Theorem \ref{T0} (2) for the sequence $C_{n,1}= B_nC_n$,
and for $H=H^{(\a)}$ and $H=u_n/2$, and  recalling \nref{test3}, we
obtain
$$
\inf_{f\in\CF(r_n)}h_n(f)\ge u_n(1+o(1)).
$$
By \nref{test2}, Corollary 2 (2) now follows.

In order to get the lower bounds, observe first that asymptotics of
$u_n$ do not depend on a sequence $C_n$ involved in \nref{balance}.
In fact, if $C_{n,0}$ is another sequence applicable to
\nref{balance}, then $C_{n,0}\sim B_nC_n,\ B_n\asymp 1$ and, under
{\bf (B3)}, we have $N(C_{n,0})\sim N(C_{n})$. Fix now a sequence
$C_n$ in \nref{balance}. It suffices to consider the case $u_n\asymp
1$, which corresponds to having $r_nC_n\sim A_n\asymp 1$. By taking
another sequence $C_{n,0}= B_nC_n,\ B_n\sim (2A_n)^{-1}$, we get
$r_nC_{n,0}\sim 1/2$. Applying Theorem \ref{T0} (1), Corollary 2 (1)
now follows. This completes the proof of Corollary 2. \endproof

\section{Sharp optimality}

\subsection{Extremal problem}
In order to describe the sharp asymptotics similar to \cite{I.93},
\cite{IS.02}, we have to consider an extremal problem on the space
of collections $v=\{v_l\}_{l\in\CL}$.

Assume that $r_n\to 0$. For $b=b_n\asymp 1, B=B_n\asymp 1$, by
following arguments similar to those in Chapter 4 of \cite{IS.02},
we arrive at
\begin{eqnarray}\label{E.2}
&&u^2_n(b,B)=\inf_{v\in V_n(b,B)}\frac 12 \sum_{l\in
\CL}v_{l}^4,\\\label{E.2a} &&V_n(b,B)=\bigg\{v :\sum_{l\in \CL}
v_{l}^2\geq n(Br_n)^2,\quad \sum_{l\in \CL} c_{l}^2 v_{l}^2\leq
nb^2\bigg\}.
\end{eqnarray}
Let $u_n(B) = u_n(1,B)$ and $\ u_n = u_n(1,1)$. From Proposition 2.8
of \cite{IS.02}, it follows that $u^2_n(b,B)$ is a convex function
in $(b^2,B^2)$ and, from rescaling arguments, it is easily seen that
$ u^2_n(b,B)=b^4u^2_n(B/b).$

By using Lagrange multipliers, the extremal collection
$v_n=\{v_{l,n}\}_{l\in\CL}$ in \nref{E.2} is of the form $
v_{l,n}^2=z_0^2(1-(c_l/C)^2)_+$, where $a_+=\max(0,a)$ for any real
number $a$, and the quantities $z_0 = z_{n,0}(b,B)>0,\ C = C_n(b,B)$
are determined by the equations
\begin{eqnarray}\label{eq.1}
&&\sum_{l\in  \CL} v_{l,n}^2=z_0^2\sum_{c_l<C}(1-(c_l/C)^2)=
n(Br_n)^2,\\
\label{eq.2} &&\sum_{l\in  \CL}c_l^2
v_{l,n}^2=z_0^2\sum_{c_l<C}c_l^2(1-(c_l/C)^2)= nb^2,
\end{eqnarray}
while the value of the extremal problem is
\begin{equation}\label{u}
u_n^2(b,B)=\frac 12 \sum_{l\in  \CL} v_{l,n}^4=\frac 12
z_0^4\sum_{c_l<C}(1-(c_l/C)^2)^2.
\end{equation}
Let
\begin{eqnarray*}
I_1&=& \sum_{l \in \CN}(1-(c_l/C)^2),\quad I_0=\sum_{l \in
\CN}(1-(c_l/C)^2)^2,\\ I_2&=& \sum_{l \in
\CN}(c_l/C)^2(1-(c_l/C)^2).
\end{eqnarray*}
It is easily seen that the equations \nref{eq.1}--\nref{u} can be
rewritten in the form
\begin{equation}\label{CB}
z_0^2I_1=n(Br_n)^2,\quad C^2z_0^2I_2=nb^2, \quad u_n^2(b,B)=\frac 12
z_0^4I_0=\frac{n^2(Br_n)^4 I_0}{2I_1^2}.
\end{equation}
Observe that $I_1=I_0+I_2\ge I_2$ and
$$
C^2=\frac{b^2I_1}{I_2B^2r_n^2}\ge b^2(Br_n)^{-2}\to\infty\quad
\text{as}\quad r_n \to 0.
$$
Under ({\bf A1}), this yields $N\to\infty$. Moreover, one has
$$
(3/4)N(C/2)\le I_1\le N(C),\quad (3/4)^2N(C/2)\le I_0\le N(C).
$$
Hence, under ({\bf B2}), these yield
\begin{equation}\label{cont}
I_1\asymp I_0\asymp N,\quad z_0^2\asymp \frac{nr_n^2}{N}, \quad
u_n^2(b,B)\asymp \frac{n^2r_n^4}{N}.
\end{equation}
Introduce the additional assumption

\medskip

({\bf C1})\quad For all $B=B_n\asymp 1$, $u_n(B)\asymp  u_n$.

\medskip

Note that, under assumption ({\bf C1}), we get
\begin{equation*}\label{cont1}
u^2_n(b,B)\sim u^2_n \quad \text{as}\,\, \ b=b_n\to 1,\ B=B_n\to 1.
\end{equation*}
(compare with Propositions 2.8 and 5.6 in \cite{IS.02}).

\subsection{Sharp asymptotics}

\begin{theorem}\label{T1} $\ $ Let $r_n\to 0$.

(1) [Lower bounds] Assume ({\bf A1})--({\bf A2}), ({\bf B1})--({\bf
B2}) and ({\bf C1}). Then
\begin{equation}\label{sh.1}
\b( \CF,r_n,\a)\ge \Phi(H^{(\a)}-u_n)+o(1),\quad\g( \CF,r_n)\ge
2\Phi(-u_n/2)+o(1),
\end{equation}
where $u_n$ is the value of the extremal problem \nref{E.2},
\nref{E.2a} for $b=B=1$.

(2) [Upper bounds] Assume ({\bf A1})--({\bf A3}) and ({\bf
B1})--({\bf B2}). Let $\liminf u_n>0$. Consider the sequence of
tests $ \psi_{n}^{H}={\1}_{\{U_n>H\}} $ based on the $U$-statistics
$$
U_n = \frac{1}{n}\sum_{1\le i<k\le n}K_n(z_i,z_k),
$$
where $z_i = (x_i,t_i),\ i=1,\ldots,n$, are the observations, with
the kernel
\begin{equation}\label{test}
K_n(z^{'},z^{''}) = x^{'}x^{''}G_n(t^{'},t^{''}),\quad
G_n(t^{'},t^{''}) = \sum_{l \in \CN}
w_{n,l}\phi_l(t^{'})\phi_l(t^{''}),
\end{equation}
where $w_{n,l}=v_{l,n}^2/u_n$ and $\{v_{l,n}\}$ is the extremal
sequence of the extremal problem \nref{E.2}, \nref{E.2a} for
$b=B=1$, or, equivalently,
$$
w_{n,l}=(1-(c_l/C)^2)_+/w_n,\quad w_n^2=\frac 12\sum_{l \in
\CN}(1-(c_l/C)^2)^2.
$$
Then, uniformly over $H=H_n\in \R$,
$$
\a(\psi_{n}^{H})\le 1-\Phi(H)+o(1),
$$
and, for any $c\in (0,1)$, uniformly over $H=H_n$ such that $u_n\ge
cH_n$,
\begin{equation}\label{sh.2}
\b(\CF,r_n,\psi_{n}^{H})\le \Phi(H-u_n)+o(1).
\end{equation}
\end{theorem}

\begin{remark}\label{R2*}
{\rm Combining  \nref{sh.1} and \nref{sh.2}, we see that the
sequence of tests $\psi_{n}^{H}$ with $H=H^{(\a)}$ is asymptotically
minimax under the Neyman-Pearson criterion, i.e.,
$$
\a(\psi_{n}^{H^{(\a)}}) \leq \a+o(1),\quad
\b(\CF,r_n,\psi_{n}^{H^{(\a)}})= \Phi(H^{(\a)}-u_n)+o(1),
$$
and the sequence of tests $\psi_{n}^{H}$ with $H=u_n/2$ is
asymptotically minimax under the total error probability criterion,
i.e.,
$$
\g( \CF,r_n,\psi_{n}^{u_n/2})= 2\Phi(-u_n/2)+o(1).
$$}
\end{remark}

\section{Tensor product Fourier basis}

Let $\Z^\infty_*\subset \Z^\infty$ consists of all sequences
$l=(l_1,\ldots,l_d,\ldots)$ with finite number $j$ such that
$l_j\not=0$, and consider the natural embedding $\Z^d\subset
\Z^\infty_*: (l_1,\ldots,l_d)\to (l_1,\ldots,l_d, 0,\ldots)$. Let
$\CL$ be an infinite subset of $\Z^\infty_*$.

Consider the tensor product Fourier basis $\{\phi_l\}_{ l\in\CL}$ in
$L_2$, i.e.,
\begin{equation}
\label{eq:tpfF} \phi_l(t) = \prod_k\phi_{l_k}(t^k),\quad
t=(t^1,\ldots t^d,\ldots)\in \Delta,\quad l\in \CL,
\end{equation}
where $\phi_j(u),\ j\in\Z, \ u\in [0,1]$, is the standard Fourier
basis in $L_2([0,1])$, i.e.,
$$
\phi_0(u) = 1,\quad \phi_j(u) = \sqrt{2}\cos(2\pi ju),\quad
 \phi_{-j}(u) = \sqrt{2}\sin(2\pi ju), \quad j>0.
$$

\begin{definition} A set $\CL$ is called {\it sign-symmetric} if, for all
$l=(l_1,\ldots,l_d,\ldots)\in\CL$, one has $\e
l=(\e_1l_1,\ldots,\e_dl_d,\ldots)\in\CL$ for all $\e_j=\pm 1$.
\end{definition}

\begin{definition}
The collection $\{h_l\}_{l \in \CL}$ is called {\it sign-symmetric}
if the set $\CL$ is sign-symmetric and $h_l=h_{\e l}$  for all
$l\in\CL$ and $\e=(\e_1,\ldots,\e_d,\ldots),\ \e_j=\pm 1$.
\end{definition}

{\bf (D1)} The set $\CL$ and the collection of coefficients
$\{c_l\}_{l \in \CL}$ are sign-symmetric.

\medskip

Let us now show  that,  under assumptions {\bf (A1)} and {\bf (D1)},
assumption {\bf (A2)} holds true for the tensor product Fourier
basis \nref{eq:tpfF}. Since the set $\CN$ is sign-symmetric then,
under assumption {\bf (D1)}, this follows from the following
statement.

\begin{lemma}
\label{lemmss} Let $\CM\subset \Z^\infty_*$ be a finite
sign-symmetric set and let $\{\phi_l\}_{l \in \CL}$ be the tensor
product Fourier basis (\ref{eq:tpfF}). Then
$$
\sum_{l\in\CM}\phi_l^2(t)= \#(\CM)\quad \forall\, t\in\Delta.
$$
\end{lemma}
{\bf Proof}. Consider the presentation $\CM=\cup_u\CM_u$, where
$u\subset \N$ and $\CM_u$ consists of $l\in\CM$ such that $\#\{j:
l_j\not=0\}=m$. It suffices to check that, for all $u$,
$$
\sum_{l\in\CM_u}\phi_l^2(t)= \#(\CM_u)\quad \forall\ t\in\Delta.
$$
Clearly, this holds for $u=\emptyset$. Without loss of generality,
assume $m=\{1,\ldots,d\},\ d\in \N$. Let $\CM_u^+=\{l\in \CM_u:
l_j>0\ \forall\ j\in u\}$. Since $\CM$ is sign-symmetric, $\CM_u^+$
consists of all $\bar\e l,\ l\in \CM_u^+$,
$\bar\e=(\e_1,\ldots,\e_d),\ \e_k=\pm 1$ and
$\#(\CM_u)=2^d\#(\CM_u^+)$. It suffices then to check that, for each
$l\in \CM_u^+$,
$$
\sum_{\bar\e}\phi_{\e l}^2(t)=2^d.
$$
Consider $\e_k$, $k=1,\ldots,d$, as $iid$ Rademacher random
variables, i.e., $P(\e_k=1)=P(\e_k=-1)=1/2$. Then, by independency,
$$
\sum_{\bar\e}\phi_{\e l}^2(t)=2^d
E_{\bar\e}\prod_{k=1}^d\phi_{\e_kl_k}^2(t^k)=2^d\prod_{k=1}^dE_{\e_k}\phi_{\e_kl_k}^2(t^k)=2^d,
$$
since
$E_{\e_k}\phi_{\e_kl_k}^2(t^k)=(2\sin^2(l_kt^k)+2\cos^2(l_kt^k))/2=1$.
This completes the proof of Lemma \ref{lemmss}.
\endproof

\begin{remark}{\rm
 Note that for the tensor product Fourier basis
(\ref{eq:tpfF}), condition \nref{sup} (and, hence, assumption {\bf
(A3)}) is fulfilled if
\begin{equation}\label{sup1}
\sum_{l\in\CL}2^{J(l)}c_l^{-2} < \infty,\quad J(l)=\#\{j:
l_j\not=0\}.
\end{equation}
Indeed, we have $ \sup_{t\in\Delta}|\phi_l(t)|=2^{J(l)/2}, $ and
hence
\begin{eqnarray*}
\|f\|_{\infty}^2\le\l(\sum_{l\in\CL}|\t_l|\sup_{t\in\Delta}|\phi_l(t)|\r)^2
&\le& \l(\sum_{l\in\CL}\t_l^2c_l^2\r)
\l(\sum_{l\in\CL}2^{J(l)}c_l^{-2}\r) \\ &\le&
\sum_{l\in\CL}2^{J(l)}c_l^{-2}.
\end{eqnarray*}
}
\end{remark}

\section{Examples: rate and sharp asymptotics in various ellipsoids}

Let us first give some extra notation. For the function
$f=\sum_{l\in \CL}\t_l\phi_l\in L_{2}^{\CL}$, we set
$\|f\|_c^2=\sum_{l\in \CL}\t_l^2c_l^2$ and let $L_{2,c}^{\CL}=\{f\in
L_{2}^{\CL}: \|f\|_c<\infty\}$ be the Hilbert space with the norm
$\|\cdot\|_c$. (Clearly the ellipsoid $\CF$ is the unit ball in
$L_{2,c}^{\CL}$.)

\medskip Consider the tensor product Fourier basis (\ref{eq:tpfF}).
In all examples below, assumption {\bf (D1)} holds true. Hence, by
Lemma \ref{lemmss}, assumption {\bf (A2)} holds true. It is easily
seen that assumption {\bf (A1)} is also fulfilled in all examples
below. That the assumption {\bf (A3)} holds also true is discussed
in each example separately.

\medskip
The first two examples are versions of the classical
multidimensional Sobolev norm (see \cite{IS.05}).

\subsection{Multidimensional Sobolev norms} Let $\Delta = [0,1]^d,\ d\in\N,\
\CL=\Z^d\setminus \{0\}$, and let
\begin{equation}\label{s1}
c_l^2 = \sum_{k=1}^d|2\pi l_k|^{2\sigma},\,\, l \in \CL, \,\,
\sigma>0.
\end{equation}
Then, for $\sigma \in\N$, the norm $\|f\|_c$ corresponds to the sum
of $\sigma$-derivatives of a $1$-periodic $f$ over all variables,
i.e.,
\begin{equation}
\label{eq:snFANIS1} \|f\|_c^2=\sum_{k=1}^d\|\partial^\sigma
f/\partial t_k^\sigma\|^2,
\end{equation}
where $\|\cdot\|$ is the norm in $L_2(\Delta)$.

Assumption {\bf (A3)} is fulfilled for $\sigma>d/4$ by the so-called
Sobolev embedding theorem (see Eq. (3.2.20) of \cite{Cohen.2003}).

Let now
\begin{equation}\label{s2}
c_l^2 = \left(\sum_{k=1}^d(2\pi l_k)^2\right)^{\sigma},\,\, l \in
\CL, \,\, \sigma>0.
\end{equation}
Then, for $\sigma \in \N$, the norm $\|f\|_c$  corresponds to the sum
of all the derivatives of a $1$-periodic $f$ of order $\sigma$,
i.e.,
\begin{equation}
\label{eq:snFANIS2} \|f\|_c^2=\sum_{i_1=1}^d \ldots
\sum_{i_\sigma=1}^d\|\partial^\sigma f/\partial
t_{i_1}\ldots\partial t_{i_\sigma}\|^2.
\end{equation}

Certainly, the norms (\ref{eq:snFANIS1}) and
(\ref{eq:snFANIS2}) are equivalent for any fixed $d$ since the ratio
of coefficients in \nref{s1} and \nref{s2} is bounded and away from
$0$. Hence, assumption {\bf (A3)} is fulfilled for $\sigma>d/4$.

It was shown in \cite{IS.05} that
$$
N(C)\sim C^{d/\sigma}J_k(d,\sigma),\ k=1,2,
$$
(e.g., $k=1$ corresponds to (6.1) and   (6.2), and $k=2$ corresponds to (6.3) and (6.4)), where
$$
J_1(d,\sigma) =
\frac{\Gamma^d(1+1/2\sigma)}{\pi^d\Gamma(1+d/2\sigma)},\quad
J_2(d,\sigma) = \frac{1}{2^d\pi^{d/2}\Gamma(1+d/2)}.
$$
Using equation \nref{balance}, these yield
$$
C\asymp n^{2\sigma/(4\sigma+d)},\quad N(C)\asymp n^{2d/(4\sigma+d)}.
$$
Hence, assumption ({\bf B2}) is fulfilled while assumption ({\bf
B1}) is fulfilled for $\sigma>d/4$. Thus, we obtain the separation
rates
$$
r_n^*=n^{-2\sigma/(4\sigma+d)}.
$$
For the sharp asymptotics, it was shown that
$$
u_n^2\sim C_k(d,\sigma)n^2r_n^{4+d/\sigma},\quad k=1,2,
$$
where, for the norm (\ref{eq:snFANIS1}),
$$
C_1(d,\sigma) =
\frac{\pi^d(1+2\sigma/d)\Gamma(1+d/2\sigma)}{(1+4\sigma/d)^{1+d/2\sigma}\Gamma^d(1+1/2\sigma)},
$$
and for the norm (\ref{eq:snFANIS2}),
$$
C_2(d,\sigma) =
\frac{\pi^d(1+2\sigma/d)\Gamma(1+d/2)}{(1+4\sigma/d)^{1+d/2\sigma}\Gamma^d(3/2)}.
$$
Assumption ({\bf C1}) is thus fulfilled. Hence, we arrive at
\nref{G}.

%

\medskip

The next two examples correspond to tensor product norms in ANOVA
modeling. These spaces are capable of dealing with interactions of
all orders in a flexible way, thus vastly extending the classical
additive methodology in multivariate nonparametric regression
inference (see \cite {Hu}, \cite{Lin}).

\subsection{Tensor product Sobolev norm} Let $\Delta = [0,1]^d,\
d\in\N,\ \CL=\Z^d$, and let
\begin{equation}\label{CAnova}
c_l = \prod_{k: l_k\not= 0}|2\pi l_k|^{\sigma}, \,\, l \in \CL, \,\,
c_{0,\ldots,0}=1.
\end{equation}
For a $\sigma \in \N$, this corresponds to the following (see
\cite{Lin}). Let us consider the functional orthogonal ANOVA
expansion
\begin{equation}\label{Anova}
f(t)=\sum_{u}f_u(t_u),\quad \int_\Delta f_u(t_u)dt_k= 0\quad
\forall\ k\in u,
\end{equation}
where the sum is taken over all subsets $u=\{j_1,\ldots
j_m\}\subset\{1,\ldots,d\}$, $1\le j_1<\ldots<j_m\le d\}$ and
$t_u=\{t_{j_1},\ldots,t_{j_m}\}$, if $u=\emptyset$, then $f_u={\rm
constant}=\int_{\Delta}f(t)dt$. Then,
$$
\|f\|_c^2=\sum_{u}\|f_u\|_{c,u}^2,
$$
where $\|f_u\|_{c,u}$ is the norm of mixed $m\sigma$-derivatives of
a $1$-periodic $f_u$, i.e.,
\begin{equation}\label{AnovaN}
\|f_u\|_{c,u}=\|\partial^{m\sigma} f/\partial
t_{j_1}^{\sigma}\ldots\partial t_{j_m}^{\sigma}\|.
\end{equation}
Assumption {\bf (A3)} is fulfilled  for $\sigma>1/4$, using
appropriate embedding properties (see Chapter III of \cite{Teml}).

It was shown in \cite{IS.07} that
\begin{equation}\label{NTen}
N(C)\sim
\frac{C^{1/\sigma}\log^{d-1}(C)}{\pi^d\sigma^{d-1}\Gamma(d)}.
\end{equation}
Using equation \nref{balance}, this yields
$$
C\asymp \l(\frac{n^{2}}{\log^{d-1}(n)}\r)^{\sigma/(4\sigma+1)}.
$$
Hence, assumption ({\bf B2}) is fulfilled while assumption ({\bf
B1}) is fulfilled for $\sigma>1/4$. Thus, we obtain the separation
rates
$$
r_n^*= \l(\frac{\log^{d-1}(n)}{n^{2}}\r)^{\sigma/(4\sigma+1)}.
$$
For the sharp asymptotics, it was shown that
\begin{equation}\label{UTen}
u_n^2\sim
\frac{C(d,\sigma)n^2r_n^{4+1/\sigma}}{\log^{d-1}(r_n^{-1})},
\end{equation}
where
\begin{equation}\label{CD}
C(d,\sigma) =
\frac{2b(\sigma)\Gamma(d)(\pi\sigma)^d}{(1+4\sigma)^{b(\sigma)}},\quad
b(\sigma) = \frac{2\sigma+1}{2\sigma}.
\end{equation}
Assumption ({\bf C1}) is thus fulfilled. Hence, we arrive at
\nref{G}.

\subsection{ANOVA subspaces} Let $\Delta = [0,1]^d,\ d\in\N$. Taking
$m\in \{0,1,\ldots,d\}$, let $\CL_{m}^{d}$ be the set that consists
of $l\in \Z^d$ such that $\#\{k: l_k\not= 0\}=m$, and
$\CL^{d,m}=\bigoplus_{j=0}^m\CL_{j}^{d}$. Under \nref{Anova}, the
spaces $L_2^{\CL^{d}_{m}}$ and $L_2^{\CL^{d,m}}$ consist of the
functions
$$
f(t)=\sum_{u: \#(u)=m}f_u(t_u),\quad f(t)=\sum_{u: \#(u)\le
m}f_u(t_u),
$$
respectively, i.e., they consist of sums of functions of $m$
variables or no more than $m$ variables. If $m=0$, this corresponds
to the constant function while the case $m=1$ corresponds to
functions with an additive structure. Take $c_l$ according to
\nref{CAnova}. Then, we obtain,
$$
\|f\|_c^2=\sum_{u: \#(u)=m}\|f_u\|_{c,u}^2,\quad \|f\|_c^2=\sum_{u:
\#(u)\le m}\|f_u\|_{c,u}^2,
$$
respectively, where, for $\sigma \in \N$, the norm $\|f_u\|_{c,u}$
of a $1$-periodic $f_u$ is determined by \nref{AnovaN} (see
\cite{Lin}). Assumption {\bf (A3)} is fulfilled for $\sigma>1/4$,
since the spaces presented here are subspaces of the tensor product
Sobolev spaces discussed in Section 6.3.

Take $c_l$ according to \nref{CAnova}. Denote by $N_d(C)$ the
function $N(C)$ for the tensor product Sobolev norms, by
$N_{d,m}(C)$  the function $N(C)$ for $\CL=\CL^{d,m}$, and by
$N_{m}^{d}(C)$ the function $N(C)$ for $\CL=\CL_{m}^{d}$. Observe
that
$$
N_{m}^{d}(C)=\binom{d}{m}N_m^m(C),\quad N_{d,m}(C)=\sum_{j=0}^m
\binom{d}{j} N_j^d(C).
$$

Set $M = \binom{d}{m}$ and note that $M\geq 1$ for $0 \leq m \leq
d$. It was shown in \cite{IS.07} that, as $C \rightarrow \infty$,
\begin{equation}
\label{eq:ANOVAasyFANIS1} N_{d,m}(C) \sim M N_m^m(C) \sim M N_m(C)
\sim \frac{M C^{1/\sigma} \log^{m-1}(C)}{\pi^m \sigma^{m-1}
\Gamma(m)},
\end{equation}
the last relation follows from \nref{NTen}. For both the cases
$\CL_{m}^{d}$ and $\CL^{d,m}$, using \nref{balance}, we have
$$
C\asymp
\l(\frac{\tilde{n}^{2}}{\log^{m-1}(\tilde{n})}\r)^{\sigma/(4\sigma+1)},
\quad \tilde{n} \eq n /\sqrt{M}.
$$
Hence, assumption ({\bf B2}) is fulfilled while assumption ({\bf
B1}) is fulfilled for $\sigma>1/4$. Thus, we obtain the separation
rates
$$
r_n^*= \l(\frac{\log^{m-1}(\tilde n)}{{\tilde
n}^{2}}\r)^{\sigma/(4\sigma+1)}.
$$
Let $u_{n,d}$ be the quantities that determine the sharp asymptotics
for the tensor product Sobolev norms with sharp asymptotics
\nref{UTen}. Using \nref{eq:ANOVAasyFANIS1},  we obtain, for both
cases, the sharp asymptotics
\begin{equation}\label{loss}
u_n^2\sim \frac{u_{n,m}^2}{M}\sim
\frac{C(m,\sigma)n^2r_n^{4+1/\sigma}}{M\log^{m-1}(r_n^{-1})},
\end{equation}
where the constant $C(m,\sigma)$ is defined by \nref{CD}. (Note that
\nref{loss} corresponds, in the case $m<d$, to some loss of
efficiency compared to \nref{UTen}, since the sample size $n$ is now
reduced by the factor $M^{-1/2} > 1$.) Assumption ({\bf C1}) is thus
fulfilled. Hence, we arrive at \nref{G}.

%

\medskip
The next example corresponds to classical multivariable analytic
functions on the complex strip (see \cite{K}, \cite{LS}).

\subsection{Multivariable analytic functions on the complex strip} Let
$\Delta = [0,1]^d,\ d\in\N,\ \CL=\Z^d$ and, for $\kappa>0$, let
$$
c^2_l = \prod_{k=1}^d\cosh(2\pi\kappa l_k), \,\, l \in \CL.
$$
This corresponds to analytic functions $f$ that provide periodic
extensions to the complex $d$-dimensional strip $(t_1+iu_1,\ldots,
t_d+iu_d),\ |u_k|\le \kappa$ (i.e., of size $2 \kappa$), and
$$
\|f\|_c^2=2^{-d}\sum_{\bar\e}\|f(\cdot+\e_k\kappa)\|^2.
$$
This case is closely related to the case
$$
c^2_l = \exp\l(2\pi\kappa\sum_{k=1}^d  |l_k|\r),\,\, l \in \CL
$$
(see \cite{LS}). Using $e^{|x|}/2\le \cosh(x)\le e^{|x|}$, condition
\nref{sup} is fulfilled for any $\kappa>0$ (by Remark 5.1), since
$$
\sum_{l\in\CL}2^{J(l)}c_l^{-2}\le
2^d\sum_{l\in\CL}c_l^{-2}\l(1+2\sum_{k=1}^\infty\exp(2\pi\kappa
k)\r)^d<\infty.
$$
Thus,  assumption {\bf (A3)} is fulfilled.

It was shown in \cite{IS.07} that
$$
N(C)\sim  \frac{2^d \log^d(C)}{(\pi \kappa)^d\Gamma(d+1)}.
$$
Using equation \nref{balance}, this yields
$$
C\asymp \frac{n^{1/2}}{(\log(n))^{d/4}}.
$$
Hence, assumptions ({\bf B1}), ({\bf B2}) are fulfilled; moreover
$N(C)$ is a slowly varying function, i.e., assumption ({\bf B3}) is
also fulfilled. Thus, we get the separation rates
$$
r_n^*=\frac{(\log(n))^{d/4}}{n^{1/2}},
$$
and the sharp asymptotics
$$
u_n^2\sim \frac{(\pi\kappa)^d\Gamma(d+1)n^2r_n^4}{2\log^d(n)}.
$$
Assumption ({\bf C1}) is thus fulfilled. Hence, we arrive at
\nref{G}.

%

\medskip
The last example corresponds to an infinitely dimensional extension
of the ANOVA decomposition, that was first suggested to lift the
curse of dimensionality in high-dimensional numerical integration
(see \cite{KS.05}, \cite{SW}, \cite{W.06}).

\subsection{Sloan-Wo$\rm\acute{z}$niakowski norm}

Let $\Delta = [0,1]^\infty,\ \CL=\Z^\infty_*$. Taking $\sigma>0, \
s>0$, let
$$
c_l = \prod_{j\in\N\,:\, l_j\not=0}j^{s}|2\pi l_j|^{\sigma},\quad
l\in \CL,\,\, s>0,\,\, \sigma>0,\,\, c_{0,\ldots,0,\ldots}=1.
$$
This corresponds to an infinite tensor product of weighed Hilbert
spaces. Under an infinite-dimensional ANOVA expansion,
$$
f(t)=\sum_{u}f_u(t_u),\quad \int_\Delta f_u(t_u)dt_k= 0\quad
\forall\ k\in u,
$$
where the sum is taken over all finite subsets $u\subset\N$, we
obtain
$$
\|f\|_c^2=\sum_{u}\gamma(u)\|f_u\|_{c,u}^2,\quad
\gamma(u)=\prod_{k\in u}k^{2s},
$$
and, for $\sigma \in \N$, the norm $\|f_u\|_{c,u}^2$ of a
$1$-periodic $f_u$ is determined by \nref{AnovaN} (see \cite{IS.06}
and compare with \cite{KS.05}, \cite{SW}, \cite{W.06}).

Contrary to the previous examples, we are not aware of any embedding
theorems for spaces of the Sloan-Wo$\rm\acute{z}$niakowski type, and
hence we cannot verify Assumption ({\bf A3}) under minimal
smoothness conditions (like $\sigma^* \eq \min(\sigma,s)>1/4$).
However, condition \nref{sup}, which leads to the Assumption ({\bf
A3}), is fulfilled for $\sigma^*>1/2$. Indeed, let $(x_{k,j}),\ k\in
\Z, 1\le j\le d $, be a matrix. Applying the formula
$$
\sum_{\bar l\in \Z^d}\prod_{j=1}^d
x_{_{l_j,j}}=\prod_{j=1}^d\sum_{l\in \Z}x_{k,j},\quad \bar
l=\{l_1,\ldots,l_d)\in \Z^d,
$$
to the matrix entries
$$
x_{k,j}=\begin{cases} 1, & k=0,\\
2j^{-2s}|2\pi k|^{-2\sigma}, & k\not=0,
\end{cases}
$$ and letting $d\to\infty$, we get, for $\sigma>1/2$ and
$s>1/2$,
\begin{eqnarray*}
\sum_{l\in\CL}2^{J(l)}c_l^{-2}&=&\sum_{l\in\CL}\prod_{j\in\N\,:\,
l_j\not=0}2j^{-2s}|2\pi
l_j|^{-2\sigma}\\
&=&\prod_{j\in\N}\l(1+2j^{-2s}\sum_{k\in \breve{\Z}}|2\pi
k|^{-2\sigma}\r)<\infty; \quad \breve{\Z}=\Z\setminus \{0\}.
\end{eqnarray*}
Thus, by Remark 5.1, assumption {\bf (A3)} is fulfilled for
$\sigma^*>1/2$.

For simplicity, we consider below only the case $\sigma\not=s$. It was shown in \cite{IS.06} that if $0<\sigma<s$ then
$$
N(C)\sim A_1C^{1/\sigma}\exp(A_2(\log C)^{\sigma/(\sigma+s)})(\log
C)^{-A_2}, $$ 
and that if $0<s<\sigma$ then $$ 
N(C)\sim B_1 C^{1/s}\exp(B_2(\log C)^{1/2})(\log C)^{-B_3},$$
where $A_i$, $i=1,2$, and $B_i$, $i=1,2,3$, are positive constants
which only depend on $\sigma, s$. Recall  that
$\sigma^* \eq\min(s,\sigma)$. Then, we get the following
log-asymptotics
$$
\log(N(C))\sim \frac{\log(C)}{\sigma^*},
$$
which correspond to the Sobolev norms for $d=1$ and
$\sigma=\sigma^*$.

It also follows that assumption ({\bf B2}) is fulfilled while
assumption ({\bf B1}) is fulfilled for $\sigma^*>1/4$. The
separation rates are of the following form. If $0<\sigma<s$, then
$$
r_n^* \asymp
n^{-2\sigma/(4\sigma+1)}\exp\left(C_1(\log
(n))^{\sigma/(s+\sigma)}\right)(\log (n))^{-C_2},
$$
and if $0<s<\sigma$, then
$$
r_n^* \asymp n^{-2s/(4s+1)}\exp\left(D_1\sqrt{\log
(n)}\right)(\log (n))^{-D_2}.
$$
These yield the following log-asymptotics
$$
\log(r_n^*)\sim - \frac{2\sigma^*\log(n)}{4\sigma^*+1}.
$$
The sharp asymptotics are of the following form. If $0<\sigma<s$,
then
$$
u_n^2\sim  C_3 n^2r_n^{4+1/\sigma}\exp\left(-C_4(\log
r_n^{-1})^{\sigma/(s+\sigma)}\right)(\log r_n^{-1})^{C_5}.
$$
If $0<s<\sigma$, then
$$
u_n^2\sim D_3 n^{2} r_n^{4+1/s}\exp\left(-D_4\sqrt{\log
r_n^{-1}}\right)(\log r_n^{-1})^{3/4},
$$
where $C_i$, $i=1,\ldots,5$, and $D_i$, $i=1,\ldots,4$, are positive
constants which only depend on $\sigma, s$. Thus, assumption ({\bf
C1}) is fulfilled. Hence, we arrive at \nref{G}.

\section{Some General Remarks}

In this section, we discuss how the main results, established in
Theorems \ref{T0} and \ref{T1} (and, hence, Corollaries 1 and 2) can be extended to
more general settings, involving non-uniform design schemes and
unknown variances. Some remarks about adaptivity issues are also
presented. We also present other, than the Fourier basis and its
tensor product version, examples of basis functions that satisfy
assumption ({\bf {A2}}), and reveal how assumption ({\bf {A2}}) can
be replaced by a weaker assumption at the cost of replacing
assumption ({\bf B1}) with a slightly stronger assumption.

\subsection{General random design schemes}

The main results, established in Theorems \ref{T0} and \ref{T1}, are
evidently extended to random design points
$y=(y^1,\ldots,y^d)\in\R^d$, $d \geq 1$, with a {{\it known}}
product probability density function, $p(y)=p_1(y^1)\times \ldots
\times p_d(y^d)$, by applying the coordinates Smirnov transform,
i.e., $y\to F(y)=(F_1(y^1),\ldots,F_d(y^d)) \in \Delta = [0,1]^d$,
where $F_k$ is the cumulative distribution function corresponding to
the probability density function $p_k$. Indeed, consider the
goodness-of-fit testing problem for testing the null hypothesis
$H_0: f=0$ against the alternative $H_1: f \in \CF_P:~\|f\|_{2,P}
\geq r_n$, where $\CF_P$ consists of functions defined on $\R^d$ and
which have the form $g(y)=f(F(y))$, $y \in \R^d$, with $g \in \CF$
and $\|f\|_{2,P} = (\int_{\R^d}f^2(y)p(y)d(y))^{1/2}$; note that, in
this case, $\|f\|_{2,P} =\|g\|$. The corresponding test statistics
are now based on the kernels \nref{test1} and \nref{test} with
$t=(t^1,\ldots,t^d)$ replaced by $F(y)=(F_1(y^1),\ldots,F_d(y^d))$
(compare with \cite{I.93a}).

We conjecture that the main results, established in Theorems
\ref{T0} and \ref{T1}, can be also extended, subject to some
additional constraints similar to \cite{I.93a}, to {{\it unknown}}
product probability density functions by replacing
$F(y)=(F_1(y^1),\ldots,F_d(y^d))$ with $
F_n(y)=(F_{n,1}(y^1),\ldots,F_{n,d}(y^d))$ in the appropriate test
statistics, where $F_{n,k}$ is the empirical distribution function
corresponding to $F_k$ for the design points
$y_{1}^k,\ldots,y_{n}^k$; this development is, however, outside the
scope of this paper.

\subsection{Unknown variance}

The results obtained in Theorems \ref{T0} and \ref{T1} are
evidently true when $\xi_i \iid \CN(0,1)$ is replaced by $\xi_i \iid
\CN(0,\tau^2)$, where $\tau^2$ is a {\em known} variance with $0 <
\tau^2 < \infty$, by multiplying $u_n$ by the factor $\tau^{-2}$ and
multiplying $r^*_n$ by the factor $\tau$, for the lower bounds, and
by multiplying the kernels \nref{test1} and \nref{test} by the
factor $\tau^{-2}$, for the upper bounds.

For an \textit{unknown} variance $\tau^2$ with $0 < \beta_1 \leq
\tau^2 \leq \beta_2 < \infty$, we replace the multiplicative factor
$\tau^{-2}$ appeared in the kernels \nref{test1} and \nref{test} by
$\tau_n^{-2}$, where $ \tau_n^2 = \sum_{i=1}^nx_i^2. $ It is easily
seen that
$$
E_{n,f}\tau_n^2=\tau^2+\|f\|^2,\quad \Var_{n,f}\tau_n^2 =
\frac{1}{n}(\|f\|^4_4-\|f\|^4 + 4 \tau^2 \|f\|^2 + 2 \tau^4)=o(1),
$$
the latter being true from assumption {\bf (A3)}. These yield
$\tau_n^2 \sim (\tau^2+\|f\|^2)$, in $P_{n,f}$-probability, which
makes possible to repeat all the arguments presented in Appendix 2
(observe that, in Appendix 2, $\|f\|^2=o(1)$ for ``least favorable''
alternative functions $f \in \CF$).

The above observations indicated that the main results established
in Theorems \ref{T0} and \ref{T1}  still remain true when the
variance $\tau^2$ is either known or, when unknown, is replaced by
an appropriate estimator as the one considered above.

\subsection{Adaptivity}\label{adapt}

Typically, the smoothness  parameter ($\sigma$ for  Sobolev norms,
$\kappa$ for analytic function, $\min(\sigma,s)$ for
Sloan-Wo$\rm\acute{z}$niakowski norms) is {\em unknown}. This leads
to the so-called problem of \textit{adaptivity}: one has to
construct a test procedure that provides the best minimax efficiency
(separation rates or sharp asymptotics) for a wide range of values
of the unknown smoothness parameter. This problem was first studied
in \cite{Spok}, and further developed in Chapter 7 of \cite{IS.02},
for the 1-variable Gaussian white noise model. The idea is to use
the Bonferroni procedure, i.e., to combine a collection of tests for
a suitable grid in a region of the unknown smoothness parameter. It
was shown in \cite{IS.02} and \cite{Spok} that this procedure
provides an asymptotically minimax adaptive testing with a small
loss (one gets an additional (but unavoidable) $\log\log(\e^{-1})$
factor in the separation rates). We conjecture that these ideas of
adaptivity could be also developed for the multivariate
nonparametric regression models considered in this paper but the
exact details should be carefully addressed; this development is,
however, outside the scope of this paper.

\subsection{Other examples of basis functions satisfying Assumption
({\bf {A2}})}

\medskip
\noindent (a) ({\em Haar basis}): Let $\phi_{jk}(t)$,
$j=0,1,\ldots$,\, $k=1,\ldots,2^j$, $t\in [0,1]$, be the standard
Haar orthonormal system on $[0,1]$ (see, e.g., Chapter 7 of
\cite{W.94}), where $j$ is the scale parameter and $k$ is the shift
parameter. Note that, in this case, $\sum_k\phi_{jk}^2(t)=2^j$, for
each resolution $j$. Consider now the tensor product version of the
Haar basis on $\Delta = [0,1]^d$, $d \geq 1$, and consider
coefficients $c_l=c_j,\ l=((j_1,k_1),\ldots,(j_d,k_d))$, that only
depend on the scale parameter $j=(j_1,\ldots,j_d)$ and not on the
shift parameter $k=(k_1,\ldots,k_d)$. Hence, by working along the
lines of Section 5, it follows that the tensor product Haar basis
functions on $\Delta$ satisfy Assumption {\bf (A2)}.

\medskip
\noindent(b) ({\em Walsh basis}): Let $\phi_j(t)$, $j=0,1,\ldots$,\,
$t \in [0,1]$, be the Walsh basis functions system on [0,1]; the
Walsh basis functions take actually sums and differences of the Haar
basis functions to obtain a complete orthonormal system (see, e.g.,
Chapter 7 of \cite{W.94}). Note that, in this case, $|\phi_j(x)|=1$,
for each $j$. Consider now the tensor product version of the Walsh
basis functions on $\Delta = [0,1]^d$, $d \geq 1$. Hence, it follows
immediately that the tensor product Walsh basis functions on
$\Delta$ satisfy Assumption {\bf (A2)}.

\medskip
\noindent (c) ({\em Orthonomal basis on a compact connected
Riemannian manifold without boundary}): Let S be a compact connected
Riemannian manifold without boundary and consider the orthonormal
system of eigenfunctions $\phi_{jk}(x)$, $x\in S$, associated with
the Laplacian (Laplace-Beltrami operator) on $S$, for different
eigenvalues $\lambda_j$, $\la_1 < \la_2 <\ldots$ with $\lambda_j
\rightarrow \infty$ as $j \rightarrow \infty$ (see, e.g.,
\cite{Gine}). For each $j=1,2,\ldots$, they satisfy the relation
$\sum_{k=1}^{k_j}(\phi_{j,k}^2(x)-\mu^{-1}(S))=0$, where $k_j <
\infty$ is the (algebraic) multiplicity of the eigenvalue
$\lambda_j$ and $\mu$ is the invariant measure on $S$ (see, e.g.,
formula (3.18), p. 127 of \cite{Efr}, or the last line of p. 1256 of
\cite{Gine}). The above relation is a natural and deep extension of
the classical relation $\sin^2(x)+\cos^2(x)=1$ for the 1-dimensional
circle. Similar to (a), consider now coefficients $c_{(j,k)}=c_j$ or
corresponding coefficients $c_l=c_j$ for the tensor product basis
functions on $S^d$, $d \geq 1$. Hence, by working along the lines of
Section 5, it follows that the tensor product basis functions on
$S^d$ satisfy Assumption {\bf (A2)}. Therefore, our general
framework could be a platform to derive analogous statements to the
ones given in Theorems \ref{T0} and \ref{T1} for minimax
goodness-of-fit testing in nonparametric regression problems on
compact connected Riemannian manifolds without boundary, $S$, or
their products, $S^d$, but the details in the derivation of these
statements should be carefully addressed; this development is,
however, outside the scope of this paper.

\subsection{Replacing assumption {\bf (A2)} by a weaker assumption}

Assumption {\bf (A2)} can be replaced by the weaker assumption
$$
{(\bf{A2a})} \qquad \sup_{t \in \Delta} \sum_{l\in \CN(C)}
\phi_l^2(t)=O(N(C)) \quad \text{as} \quad C \rightarrow \infty,
$$
(it covers the cosines orthonormal system, compactly supported
(other than the Haar basis) orthonormal wavelet systems, as well as
their tensor product versions) by replacing assumption ({\bf {B1}})
with the slightly stronger assumption
$$
({\bf {B1a}}) \qquad  N=o(n^{2/3}).
$$
Indeed, the only difference in the proofs of Theorems \ref{T0} and \ref{T1} is in the relation \nref{X}. In particular, one can
use the Cauchy-Schwarz inequality which yields an additional factor
$N$, and this is compensated by assumption ({\bf {B1a}}).

\section{Appendix 1: proof of lower bounds}

Let us start with some extra notation. Recall first that
$X_n=\{x_1,\ldots,x_n\}$, $T_n=\{t_1,\ldots,t_n\}$, $Z_n=(X_n,T_n)$,
and $z_i=(x_i,t_i)$, and that $P_{n,f}$ is the probability measure
that corresponds to $Z_n$ whereas $E_{n,f}$ is the expectation over
this probability measure. Denote also by $\Var_{n,f}$ the
corresponding variance. Let $P_{n,T}$ be the probability measure
that corresponds to $T_n$ and $P_{n,f}^T$ be the conditional
probability measure with respect to $T_n$. We denote by $E_{n,T}$
and $E_{n,f}^T$ the expectations over these probability measures,
whereas $\Var_{n,T}$, $\Var_{n,f}^T$ are the corresponding
variances. (Clearly, $E_{n,f}(\cdot)=E_{n,T} E_{n,f}^T(\cdot)$.)
Also, for the function $f=\sum_{l}\theta_l\phi_l$, we denote the
measure $P_{n,f}$ by $P_{n,\theta}$, with analogous notation for the
expectations, conditional expectations and variances. Let also
$E^{T,\xi}_{n}$ and $\Var^{T,\xi}_{n}$ be the expectation and
variance of the conditional probability measure with respect to
$\Xi_n =\{\xi_1,\ldots,\xi_n\}$, where $\xi_i\iid \CN(0,1)$.
Certainly, $P_{n,\xi}=P_{n,0}$.

\subsection{Lower bounds for Theorem \protect{\ref{T1}}}

\subsubsection{Priors} We use the constructions similar to
\cite{E.90a} and follow, but with necessary modifications,
techniques from \cite{I.93}--\cite{IS.02}. It suffices to consider
the case
\begin{equation}\label{A0}
u_n^2\asymp  1.
\end{equation}
Take $\delta\in (0,1)$, let $a_{l,n}=v_{l,n}(b,B)$ be the extremal
collection for the extremal problem \nref{E.2}, \nref{E.2a} with
$b=1-\delta,B=1+\delta$,
and let $A=A_n$ be the diagonal matrix with diagonal elements
$a_l=a_{l,n},\ l\in\CN$.

Under \nref{A0}, using ({\bf C1}), \nref{cont}, we have
\begin{equation}\label{A}
u^2_n(b,B)=\frac{1}{2}\sum_{l\in \CN}a^4_{l,n}\asymp 1,\quad D_n =
N\max_{j\in \CN}a^4_{j,n}\sim z_0^4N\asymp 1 .
\end{equation}
Let $v=\sqrt{n}\t$ and let $\pi_n(dv)$ be the Gaussian prior
$\CN(0,A^2)$ on the parametric space consisting of
$\{v_l\}_{l\in\CL}=\sqrt{n}\{\t_l\}_{l\in\CL}$, i.e., $v_l$ are
independent in $l$ and, for each $l$, $v_l \sim\CN(0,a^2_{l})$ for
$c_l<C$ and $v_l=0$ for $c_l\ge C$, in $\pi_n$-probability.

Note that, in the sequence space of the ``generalized'' Fourier
coefficients $\t=\{\t_l\}_{l\in\CL}$ with respect to the orthonormal
system $\{\phi_l\}_{ l\in\CL}$, the null hypothesis (\ref{eq:nullF})
(recall that $f_0 = 0$) corresponds to $H_0:
 \t=0$ and, assuming $f \in \CF$, the alternative hypothesis (\ref{eq:altF}) corresponds to
\begin{equation}\label{alt}
H_1 :\quad \sum_{l\in\CL}c_l^2\t_l^2\le 1,\quad
\sum_{l\in\CL}\t_l^2\ge r_n^2.
\end{equation}

Let $V_n = V_n(1,1)$ be the set determined by \nref{E.2a} with
$B=b=1$; this corresponds to the alternative set \nref{alt}.

\begin{lemma}\label{LB}
For any $\delta\in (0,1)$, one has $ \pi_n(V_n)=1+o(1). $
\end{lemma}
{\bf Proof of Lemma \ref{LB}}. It follows from evaluations of
$\pi_n$-expectations and variances of the random variables $
\Cl_1(v)=\sum_{l\in \CN}v_l^2$ and $\Cl_2=\sum_{l\in
\CN}c_l^2v_l^2$, and by using the Chebyshev inequality (compare with
similar evaluations in \cite{I.93}, \cite{IK.07}, \cite{IS.02}).
\endproof

\medskip

Let $\b(P_{n,0},P_{\pi_n},\a)$  be the minimal type II error
probability for a given level $\a\in (0,1)$ and
$\g(P_{n,0},P_{\pi_n})$ be the minimal total error probability for
testing the simple null hypothesis $H_0:P=P_{n,0}$ against the
simple Bayesian alternative $H_0:P=P_{\pi_n}$ for the mixture
$P_{\pi_n}(A)=\int P_{n,n^{-1/2}v}(A)\,\pi_n(dv)$. By Lemma \ref{LB}
and using Proposition 2.11 in  \cite{IS.02}, we have
$$
\b(\CF,r_n,\a)\ge\b(P_{n,0},P_{\pi_n},\a)+o(1),\quad
\g(\CF,r_n)\ge\g(P_{n,0},P_{\pi_n})+o(1).
$$
Hence, it suffices to show that
\begin{equation}\label{B}
\b(P_{n,0},P_{\pi_n},\a)\ge \Phi(H^{(\a)}-u_n)+o(1),\quad
\g(P_{n,0},P_{\pi_n})\ge 2\Phi(-u_n/2)+o(1).
\end{equation}
In order to obtain \nref{B}, it suffices to verify that, in
$P_{n,0}$-probability,
\begin{equation}\label{Ba}
\log(dP_{\pi_n}/dP_{n,0})=-u_n^2/2+u_n\zeta_n+\eta_n,\quad \eta_n\to
0,\quad \zeta_n\to\zeta\sim \CN(0,1)
\end{equation}
(see \cite{IS.02}, Section 4.3.1, formula (4.72)).

\subsubsection{Likelihood ratio and correlation matrix}
For $f(t)=\sum_{l\in\CN}\t_l\phi_l(t)$, the likelihood ratio is of
the form
$$
\frac{dP_{n,\t}}{dP_{n,0}}=\frac{dP^T_{n,\t}}{dP^T_{n,0}}=\exp\bigg(-\frac
12 v'Rv+\langle w,v \rangle_s \bigg),\quad \t=\{\t_l\}_{
l\in\CN},\quad v = \sqrt{n}\t,
$$
where $w=\{w_l \}_{l \in \CN}$, $w_l=w_{l,n} =
\frac{1}{\sqrt{n}}\sum_{i=1}^n x_i\phi_l(t_i)$, and $R$ is the
correlation matrix
$$
R=R_n=\{r_{jl}\}_{j,l\in \CN},\quad r_{jl} = \frac
1n\sum_{i=1}^n\phi_j(t_i)\phi_l(t_i);
$$
here, and in Section 9.1.3, $\langle \cdot,\cdot \rangle_s$ denotes
the inner product in the sequence space.

\medskip

Let $\tr(\cdot)$ be the trace of a square matrix.

\begin{lemma}\label{LR}

(1) The matrix $R$ is symmetric and positively semi-defined.
Moreover, $E_{n,T}R=I_N$, where $I_N=\{\delta_{jl}\}_{j,l\in\CN}$ is
the unit $N\times N$ matrix.

(2) Under \nref{Ph} and ({\bf B1}), one has
\begin{eqnarray}\label{R1}
 E_{n,T}\tr(R^2)&\sim& N, \\
\label{R2} E_{n,T}\tr((R-I_N)^2)&=&o(N),\\
\label{R3} E_{n,T}\tr(R^4)&\sim& N.
\end{eqnarray}
\end{lemma}
{\bf Proof of Lemma \ref{LR}}. First, we prove statement (1). For
any $\tilde x=\{\tilde x_j\}_{j\in\CN},\ \tilde x_j\in\R$, one has
$$
\sum_{j,l\in \CN}\tilde x_j \tilde x_l r_{jl}=\frac
1n\sum_{i=1}^n\left(\sum_{j\in \CN}\tilde x_j
\phi_j(t_i)\right)^2\ge 0.
$$
Since $\{\phi_l\}_{l\in \CN}$ is an orthonormal system,
$$
E_{n,T}r_{jl}=\int_\Delta\phi_j(t)\phi_l(t)dt=\delta_{jl}.
$$
Thus, statement (1) follows.

Now, we prove statement (2). Analogously, we have, using ({\bf A2}),
({\bf B1}),
\begin{eqnarray*}\nonumber
E_{n,T}(r_{jl}-\delta_{jl})^2&=&\Var_{n,T}r_{jl}=\frac
1n\left(\int_\Delta\phi_j^2(t)\phi_l^2(t)dt -\delta_{jl}^2\right)\\
\label{r2} &=&\frac
1n\int_\Delta\phi_j^2(t)\phi_l^2(t)dt-\frac{1}{n}\delta_{jl}
\end{eqnarray*}
and
\begin{eqnarray*}\nonumber E_{n,T}\tr((R-I_N)^2)&=&\sum_{j,l\in
\CN}E_{n,T}(r_{jl}-\delta_{jl})^2\le \frac
1n\int_\Delta\sum_{j,l\in \CN}\phi_j^2(t)\phi_l^2(t)dt\\
\label{tr2} &=&\frac 1n\int_\Delta\left(\sum_{j\in
\CN}\phi_j^2(t)\right)^2dt=\frac{N^2}{n}=o( N),
\end{eqnarray*}
which yields \nref{R2}. We obtain \nref{R1} from \nref{R2} since $
\tr(R^2)=\tr((R-I_N)^2)+\tr(I_N). $

Let us now evaluate $E_{n,T}\tr(R^4)$. Let $ R^2 =
\{b_{jl}\}_{j,l\in\CN},$
$$
b_{jl}=\sum_{s\in\CN}r_{js}r_{sl}=\frac{1}{n^2}\sum_{s\in\CN}\sum_{\a,\b=1}^n
\phi_j(t_\a)\phi_s(t_\a)\phi_s(t_\b)\phi_l(t_\b).
$$
We have
\begin{eqnarray*}
&&\tr(R^4)=\sum_{j,l\in\CN}b_{jl}^2\\
&&=\frac{1}{n^4}\sum_{l,j,s,r\in\CN}\sum_{\a,\b,\g,\de=1}^n
\phi_j(t_\a)\phi_s(t_\a)\phi_s(t_\b)\phi_l(t_\b)\phi_j(t_\g)\phi_r(t_\g)\phi_r(t_\de)
\phi_l(t_\de).
\end{eqnarray*}
Observe that
\begin{eqnarray*}
\sum_{\a,\b,\g,\de=1}^n&\!\!\!\!\!\!\!\!\!\!\!\!\!\!\!\!\!\!&E_{n,T}\l\{
\phi_j(t_\a)\phi_s(t_\a)\phi_s(t_\b)\phi_l(t_\b)\phi_j(t_\g)\phi_r(t_\g)\phi_r(t_\de)
\phi_l(t_\de)\r\}\\ &:=&S_4+S_3+S_2+S_1,
\end{eqnarray*}
where $S_4$--$S_1$ correspond to the sums (we omit indexes $j,l,r,s$
in notation of $S_1$--$S_4$)
\begin{eqnarray*}
S_4&=&24\sum_{1\le\a<\b<\g<\de\le n},\\
S_3&=&6\l(\sum_{1\le\a=\b<\g<\de\le n}+\sum_{1\le\a<\b=\g<\de\le n}+
\sum_{1\le\a<\b<\g=\de\le n}\r),\\
S_2&=&2\l(\sum_{1\le\a=\b=\g<\de\le n}+\sum_{1\le\a<\b=\g=\de\le
n}+\sum_{1\le\a=\b<\g=\de\le n}\r),\\
S_1&=&\sum_{1\le\a=\b=\g=\de\le n}.
\end{eqnarray*}
By independence of $t_i$, and since   $\{\phi_l\}$ is an orthonormal
system, we have
\begin{eqnarray*}
S_4&=&C_4(n)\delta_{js}\delta_{sl}\delta_{jr}\delta_{rl},\\
S_3&=&C_3(n)\Big\{\delta_{jr}\delta_{lr}\int_{\Delta}\phi_j(t)\phi_s^2(t)\phi_l(t)dt+
\delta_{js}\delta_{rl}\int_{\Delta}\phi_s(t)\phi_l(t)\phi_j(t)\phi_r(t)dt\\
&+&\delta_{js}\delta_{sl}
\int_{\Delta}\phi_j(t)\phi_r^2(t)\phi_l(t)dt
\Big\},\\
S_2&=&C_2(n)\Big\{\delta_{rl}\int_{\Delta}\phi_j^2(t)\phi_s^2(t)\phi_l(t)\phi_r(t)dt+
\delta_{sj}\int_{\Delta}\phi_l^2(t)\phi_r^2(t)\phi_j(t)\phi_s(t)dt\\
&+& \l(\int_{\Delta}\phi_j(t)\phi_s^2(t)\phi_l(t)dt\r)
\l(\int_{\Delta}\phi_j(u)\phi_r^2(u)\phi_l(u)du\r) \Big\},\\
S_1&=&n\int_{\Delta}\phi_j^2(t)\phi_s^2(t)\phi_r^2(t)\phi_l^2(t)dt,
\end{eqnarray*}
where $C_4(n)\sim n^4,\ C_3(n)\asymp n^3,\ C_2(n)\asymp n^2$.
Therefore,
\begin{eqnarray*}
\frac{1}{n^4}\sum_{l,j,s,r\in\CN}S_4&=&
\frac{C_4(n)}{n^4}\sum_{l,j,s,r\in\CN}\delta_{js}\delta_{sl}\delta_{jr}\delta_{rl}=
\frac{NC_4(n)}{n^4}\sim N,\\
\frac{1}{n^4}\sum_{l,j,s,r\in\CN}S_3&=&\frac{3C_3(n)}{n^4}\sum_{j,s\in\CN}
\int_{\Delta}\phi_j^2(t)\phi_s^2(t)dt\\
&=&\frac{3C_3(n)}{n^4}\int_{\Delta}\l(\sum_{j\in\CN}
\phi_j^2(t)\r)^2dt=\frac{3N^2C_3(n)}{n^4}=O(N^2/n),\\
\frac{1}{n^4}\sum_{l,j,s,r\in\CN}S_1&=&\frac{n}{n^4}\sum_{l,j,s,r\in\CN}
\int_{\Delta}\phi_j^2(t)\phi_s^2(t)\phi_r^2(t)\phi_l^2(t)dt\\
&=&\frac{1}{n^3}\int_{\Delta}\l(\sum_{j\in\CN}
\phi_j^2(t)\r)^4dt=\frac{N^4}{n^3}.
\end{eqnarray*}
Analogously,
\begin{eqnarray*}
&&\sum_{l,j,s,r\in\CN}\delta_{rl}\int_{\Delta}\phi_j^2(t)\phi_s^2(t)\phi_l(t)\phi_r(t)dt=
\int_{\Delta}\l(\sum_{l,j,s\in\CN}\phi_j^2(t)\phi_s^2(t)\phi_l^2(t)\r)dt\\
&&=\int_{\Delta}\l(\sum_{l\in\CN}\phi_j^2(t)\r)^3dt=N^3
\end{eqnarray*}
and
\begin{eqnarray}\nonumber
&&\sum_{l,j,s,r\in\CN}\l(\int_{\Delta}\phi_j(t)\phi_s^2(t)\phi_l(t)dt\r)
\l(\int_{\Delta}\phi_j(u)\phi_r^2(u)\phi_l(u)du\r)\\
&&=
\sum_{l,j\in\CN}\l(\int_{\Delta}\phi_j(t)\l(\sum_{s\in\CN}\phi_s^2(t)\r)\phi_l(t)dt\r)
\l(\int_{\Delta}\phi_j(u)\l(\sum_{s\in\CN}\phi_r^2(u)\r)\phi_l(u)du\r)\qquad  \label{X}\\
&&=N^2\sum_{l,j\in\CN}\l(\int_{\Delta}\phi_j(t)\phi_l(t)dt\r)\nonumber
\l(\int_{\Delta}\phi_j(u)\phi_l(u)du\r)=N^2\sum_{l,j\in\CN}\delta_{jl}^2=N^3.
\end{eqnarray}
Thus,
$$
\frac{1}{n^4}\sum_{l,j,s,r\in\CN}S_2=O(N^3/n^2).
$$
Combining evaluations above and ({\bf B1}) we get \nref{R3}:
$$
\tr(R^4)\sim N(1+O(N/n+(N/n)^2+(N/n)^3))\sim N.
$$
Thus, statement (2) follows. This competes the proof of Lemma
\ref{LR}. \endproof

\subsubsection{Bayesian likelihood ratio} Let us now study the
Bayesian likelihood ratio. Direct calculation gives
\begin{equation}\label{ll}
\frac{dP_{\pi_n}}{dP_{n,0}}=E_{\pi_n}\frac{dP^T_{n,\t}}{dP^T_{n,0}}=\frac{1}{\sqrt{\det{G}}}
\exp\left(\frac 12q'G^{-1}q\right),
\end{equation}
where $q = Aw,\ G = G_n = I_N+A'RA$. Let $\tilde \tau_l\ge 0,\ l\in
\CN,$ be the eigenvalues of the symmetric positively semi-defined
matrix $D = A'RA=\{a_ja_lr_{jl}\}_{j,l\in\CN}$. Let $e_l$ be the
eigenvectors of the matrix $D$ and let $q_l=\langle q, e_l
\rangle_s$, $l \in \CL$.

We can now rewrite \nref{ll} in the form
$$
L_n = \log\left(\frac{dP_{\pi_n}}{dP_{n,0}}\right)=\frac
12\sum_{l\in \CN}\left(\frac{q_l^2}{1+\tilde \tau_l}-\log(1+\tilde
\tau_l) \right).
$$
Let $\|\tilde{A}\|_\infty = \sup_{\|x\| \leq 1}\|\tilde{A}x\|$ for a
generic matrix $\tilde{A}$. Observe that
$$ \|D\|_\infty^4=\max_{l\in
\CN}\tilde \tau_l^4\le \sum_{l\in \CN}\tilde \tau_l^4={\tr(D^4)}.
$$
Using the standard relations
$$ \tr(AC)=\tr(CA)\quad\text{and}\quad
\tr(A'BA)\le\|A\|_\infty^2\tr(B),
$$
for a symmetric positively semi-defined matrix $B$, we get the
inequalities
$$
\tr(D^2) \le\|A\|_\infty^4\tr(R^2)\quad\text{and}\quad \tr(D^4)
\le\|A\|_\infty^8\tr(R^4).
$$
By \nref{A},
$$
\|A\|_\infty^4=\max_{l\in\CN}a_l^4\le D_n/N.
$$
Jointly with \nref{R1} and \nref{R3}, the above  yields
\begin{equation*}\label{trD}
E_{n,T}(\tr(D^2))=O(1),\quad  E_{n,T}(\tr(D^4))=O(N^{-1}).
\end{equation*}
Hence,
\begin{equation*}\label{norm}
E_{n,T}\l(\max_{l\in \CN}|\tilde \tau_l|\r)=O(N^{-1/4}).
\end{equation*}
Thus, in $P_{n,T}$-probability,
\begin{equation}\label{norm2}
\|D\|_\infty=\max_{l\in\CN}|\tilde \tau_l|=o(1).
\end{equation}

Using the well-known relations
$$ (1+y)^{-1}=1-y+o(y) \quad \text{and} \quad
\log(1+y)-y+y^2/2=o(y^2), \quad \text{as} \quad y \to 0,
$$
we get, with
$P_{n,T}$-probability tending to 1, by \nref{norm2},
\begin{eqnarray}\nonumber
L_n&=&\frac 12\sum_{l\in \CN}\left(q_l^2(1-\tilde \tau_l)-\tilde
\tau_l+\tilde \tau_l^2/2 \right)+o\left( \sum_{l\in
\CN}q_l^2\tilde \tau_l\right)+ o\left( \sum_{l\in \CN}\tilde \tau_l^2\right)\\
\nonumber &=&\frac 12\left(\tr(Q)-\tr(D)-\tr(QD)+\tr(D^2)/2\right)
+o\left(\tr(QD)\right)+o\left(\tr(D^2)\right)\\ \label{l0} &=&\frac
12\l(\tr(\hat{Q})-\tr(\hat{Q}D)-\tr(D^2)/2\r)
+o\left(\tr(\hat{Q}D)\right)+o\left(\tr(D^2)\right),
\end{eqnarray}
where
$$
Q = qq'=Azz'A=\{a_{j}a_{l}z_jz_l\}_{j,l\in \CN},\ \hat{Q} =
Q-D=A(zz'-R)A.
$$

Let us now study the $P_{n,0}$-distribution of $L_n$.
\begin{lemma}\label{L1}
In $P_{n,0}$-probability,
\begin{eqnarray}\label{l1}
\tr(\hat{Q}D)&=&o(1),\\
 \label{l2}
 \tr(D^2)&=&\tr(A^4)+o(1),\\
 \label{l3}
 E_{n,0}\tr(\hat Q)&=&0,\\
 \label{l3F}
 \Var_{n,0}\tr(\hat Q)&=&2\tr(A^4)+o(1).
\end{eqnarray}
\end{lemma}
{\bf Proof of Lemma \ref{L1}}. Let
$\Phi=n^{-1/2}\{\phi_j(t_i)\}_{j\in\CN, i=1,\ldots,n}$ be an
$N\times n$-matrix, and set $\xi'=(\xi_1,\ldots,\xi_n)$. Then, in
$P_{n,0}$-probability,
$$
R=\Phi\Phi', \quad z=\Phi\xi,\quad z'z=\xi'\Phi'\Phi\xi,\
E(\xi\xi')=I_N.
$$
Observe that
$$
E^T_{n,0}zz'=\Phi\l(E^T_{n,0}\xi\xi'\r)\Phi'=\Phi\Phi'=R,
$$
which yields
\begin{equation}\label{EQ}
E^T_{n,0}(\tr(\hat{Q}))=0,\quad E^T_{n,0}(\tr(\hat{Q}D))=0.
\end{equation}
Analogously, using the formula
$$
\Var(\tr(B\xi\xi')))=2\tr(BB'),
$$
we get
$$
\Var^T_{n,0}(\tr(\hat{Q}D))=\Var^T_{n,0}\tr(A\Phi\xi\xi'\Phi'AD)=2\tr(BB'),
$$
where $ B = \Phi'A^2\Phi\Phi'A^2\Phi.$ By Lemma \ref{LR} and
\nref{A}, it is easily seen that
$$
\tr(BB')=\tr((ARA)^4)\le\|A\|_\infty^8\tr(R^4). 
$$
Using the formula
$$
\Var_{n,0}(\cdot)=\Var_T(E^T_{n,0}(\cdot))+E_T(\Var^T_{n,0}(\cdot)),
$$
we get
$$
\Var_{n,0}(\tr(\hat{Q}D))=o(1),
$$
which together  with \nref{EQ},  yields   \nref{l1}.

To obtain \nref{l2}, note that
$$
\tr(D^2)=\tr(\hat D^2)+2\tr(A^2\hat D)+\tr(A^4),\quad \hat D =
D-A^2=A(R-I_N)A,
$$
and observe that, by Lemma \ref{LR} and \nref{A},
$$
\tr(\hat D^2)\le \|A\|_\infty^4\tr((R-I_N)^2)=o(1),\quad
(\tr(A^2\hat D))^2\le \tr(A^4)\tr(\hat D^2)=o(1).
$$

Obviously, \nref{l3} follows from \nref{EQ}, and \nref{l3F} follows
from \nref{l2}, since
$$
\Var_{n,0}^T(\tr(\hat
Q))=\Var_{n,0}^T(\tr(A\Phi\xi\xi'\Phi'A))=2\tr((A\Phi\Phi'A)^2)=2\tr(D^2).
$$
This completes the proof of Lemma \ref{L1}. \endproof

\medskip

Let $ \zeta_n = \tr(\hat Q)/2u_n,\ u_n^2=\tr(A^4)/2.$ By Lemma
\ref{L1}, we rewrite \nref{l0} in the form
$$
L_n=u_n\zeta_n-u_n^2/2+\eta_n,\quad \eta_n
\,\stackrel{P_{n,0}}{\to}\, 0.
$$

\begin{lemma}\label{L2} In $P_{n,0}$-probability,
$\zeta_n\to\zeta\sim\CN(0,1)$.
\end{lemma}
{\bf{Proof of Lemma \ref{L2}}}. Let us rewrite $\tr(\hat Q)$ in the
form
\begin{eqnarray*}
\frac 12\tr(\hat Q)&=&\frac 12\tr(A\Phi(\xi\xi'-I)\Phi'A)=\frac
12\sum_{i=1}^nw_{ii}(\xi_i^2-1)+\sum_{1\le i< k\le
n}w_{ik}\xi_i\xi_j \\
&:=&A_n+B_n,
\end{eqnarray*}
where
$$
W=\{w_{ik}\}_{i,k=1}^n=\Phi'A^2\Phi,\quad w_{ik}=\frac
1n\sum_{l\in\CN}a_l^2\phi_l(t_i)\phi_l(t_k).
$$
It is easily seen that $E^{T,\xi}_{n}A_n=0,$ and by ({\bf A2}),
\nref{A},
\begin{eqnarray*}
\Var^{T,\xi}_{n}(A_n)&=&\frac 12
\sum_{i=1}^nw_{ii}^2=\frac{1}{2n^2}\sum_{i=1}^n\l(\sum_{l\in\CN}
a_l^2\phi_l^2(t_i)\r)^2\\
&\le&\frac{D_n}{2n^2N}\sum_{i=1}^n\l(\sum_{l\in\CN}
\phi_l^2(t_i)\r)^2 = \frac{D_nN}{2n}=o(1).
\end{eqnarray*}
Thus, $A_n\to 0$ in $L_2(P_{n,0})$ and in $P_{n,0}$-probability.

The item $B_n$ is degenerate $U$-statistic
$$
B_n=\frac{1}{n}\sum_{1\le i< k\le n}W_n(r_i,r_j),\quad
r_i=(\xi_i,t_i)\quad \text{are} \quad i.i.d.,
$$
$$
W_n(r^{'},r^{''})={\xi^{'}
\xi^{''}}\sum_{l\in\CN}a_l^2\phi_l(t^{'})\phi_l(t^{''}),\quad \int
W_n(r^{'}, r^{''})P(d r^{'})=0\quad \forall r^{''},
$$
where $P(dr)=\CN_{0,1}(d\xi)\times U_{\Delta}(dt)$, i.e., $\xi$ and
$t$ are independent, $\xi\sim\CN(0,1)$ and $t$ is uniformly
distributed on $\Delta$.

The statement of Lemma \ref{L2} follows from the following
proposition.

\begin{proposition}\label{GU} In
$P_{n,0}$-probability, the statistics $B_n$ are asymptotically
$\CN(0,u_n^2)$.
\end{proposition}
{\bf Proof of Proposition \ref{GU}}. Clearly, $E_{P_{n,0}}B_n=0$
and, for $r_1=(\xi_1,t_1),\ r_2=(\xi_2,t_2)$,
\begin{eqnarray*}
\Var_{P_{n,0}}(B_n)&=&\frac{n(n-1)}{2n^2}\int\int
W_n^2(r_1,r_2)P(dr_1)P(dr_2)\\
&=&\frac{n(n-1)}{2n^2}
E(\xi_1^2\xi_2^2)\int_{\Delta}\int_{\Delta}\l(\sum_{l\in\CN}a_l^2\phi_l(t_1)\phi_l(t_2)\r)^2dt_1
dt_2\\
&=&\frac{n(n-1)}{2n^2}\sum_{j,l\in\CN}
a_j^2a_l^2\int_{\Delta}\int_{\Delta}\phi_j(t_1)\phi_j(t_2)\phi_l(t_1)\phi_l(t_2)dt_1
dt_2\\
&=& \frac{n(n-1)}{2n^2}\sum_{l\in\CN}a_l^4\sim u_n^2.
\end{eqnarray*}
For $r_1=(\xi_1,t_1),\ r_2=(\xi_2,t_2),\ r_3=(\xi_3,t_3)$, let
\begin{eqnarray*}
\widetilde{G}_n(r_1,r_2)&=&\int W_n(r_1,r_3)W_n(r_2,r_3)P(dr_3),\\
G_{n,2}&=&\int\int \widetilde{G}_n^2(r_1,r_2)P(dr_1)P(dr_2),\\
W_{n,4}&=&\int\int W^4_n(r_1,r_2)P(dr_1)P(dr_2).
\end{eqnarray*}
Using the asymptotic normality of degenerate $U$-statistics
established in \cite{H}, together with Lemma 3.4 in \cite{I.94}, it
suffices to verify the conditions
\begin{eqnarray}\label{C1}
\widetilde{G}_{n,2}&=&o(1),\\ W_{n,4}&=&o(n^2). \label{C2}
\end{eqnarray}
We have
\begin{eqnarray*}
\widetilde{G}_n(r_1,r_2)&=&E_{P(d\xi_3,dt_3)}\l(
{\xi_1\xi_2\xi_3^2}\sum_{l\in\CN}a_l^2\phi_l(t_1)\phi_l(t_3)\sum_{j\in\CN}a_j^2\phi_j(t_2)\phi_j(t_3)\r)\\
&=&{\xi_1\xi_2}\sum_{j,l\in\CN}a_l^2a_j^2\phi_l(t_1)\phi_j(t_2)\int_\Delta\phi_l(t_3)\phi_j(t_3)dt_3=
{\xi_1\xi_2}\sum_{l\in\CN}a_l^4\phi_l(t_1)\phi_l(t_2),\\
G_{n,2}&=&E(\xi_1\xi_2)^2\int_{\Delta}\int_{\Delta}\l(\sum_{l\in\CN}a_l^4\phi_l(t_1)\phi_l(t_2)\r)^2dt_1
dt_2=\sum_{l\in\CN}a_l^8=O(N^{-1}),
\end{eqnarray*}
which yields \nref{C1}. Next,
\begin{eqnarray*}
W_{n,4}&=&E(\xi_1\xi_2)^4\int_{\Delta}\int_{\Delta}\l(\sum_{l\in\CN}a_l^2\phi_l(t_1)\phi_l(t_2)\r)^4dt_1
dt_2\\
&\le&9\sup_{t_1,t_2\in\Delta}\l(\sum_{l\in\CN}a_l^2\phi_l(t_1)\phi_l(t_2)\r)^2
\int_{\Delta}\int_{\Delta}\l(\sum_{l\in\CN}a_l^2\phi_l(t_1)\phi_l(t_2)\r)^2dt_1
dt_2=O(N),
\end{eqnarray*}
since by ({\bf A2}) and \nref{A}, we have
\begin{eqnarray*}
\sup_{t_1,t_2\in\Delta}\l|\sum_{l\in\CN}a_l^2\phi_l(t_1)\phi_l(t_2)\r|=\sup_{t_1\in\Delta}
\sum_{l\in\CN}a_l^2\phi_l^2(t_1)
\le\max_{l\in\CN}a_l^2\sup_{t_1\in\Delta}\sum_{l\in\CN}\phi_l^2(t_1)=O(N^{1/2}).
\end{eqnarray*}
This implies \nref{C2}. This completes the proof of Proposition
\ref{GU}. Hence, Lemma \ref{L2} follows.
\endproof

\medskip

Thus, we obtain \nref{Ba} which yields \nref{B}. Hence, Theorem
\ref{T1} (1) follows. \endproof

\subsection{Lower bounds for Theorem \protect{\ref{T0}}}

The same scheme used in the proof of the lower bounds of Theorem 2
can be also employed here.

Let $C^2r_n^2<(1-\delta),\ \delta>0$. It suffices to assume
$u_n^2=n^2r_n^4/2N=O(1)$. We take the Gaussian prior
$\pi_n=\CN(0,A^2)$ that corresponds to the matrix $A=a_nI_N$ with
$a_n^2= nr_n^2(1+\delta)/N$. Recall $\Cl_1$, $\Cl_2$ from the proof
of Lemma \ref{LB}. Analogously to the proof of Lemma \ref{LB}, we
have
\begin{eqnarray*}
E_{\pi_n}\Cl_1&=&a_n^2N=nr_n^2(1+\delta),\\
E_{\pi_n}\Cl_2&\le&C^2a_n^2N<nC^2r_n^2(1-\delta)<n,\\
\Var_{\pi_n}\Cl_1&=&2a_{n}^4N=O(1),\\
\Var_{\pi_n}\Cl_2&\le&2C^4a_{n}^4N=O(n^2/N).
\end{eqnarray*}
Since, by Chebyshev's inequality,
$\Var_{\pi_n}\Cl_k=o((E_{\pi_n}\Cl_k)^2),\ k=1,2$, these yields
$\pi_n(V_n)=1+o(1)$.

Observe that the relations \nref{A} hold true with $z_0=a_n$.
Repeating the calculations in the proof of the lower bounds of
Theorem \ref{T1}, we arrive at \nref{B} with
$u_n^2=Na_n^4/2=n^2r_n^4/2N(1+\delta)^2$. Since $\delta>0$ can be
taken arbitrary small, this yields Theorem \ref{T0} (1).
\endproof

\section{Appendix 2: proof of upper bounds}

\subsection{Upper bounds for Theorem \protect{\ref{T1}}}

We consider the test sequence $\psi_{n}^{H}=\1_{\{U_n>H\}}$ based on
the $U$-statistics $U_n$ with the kernel $K_n(z_1,z_2)$ of the form
\nref{test}.

\subsubsection{Type I error} Observe that
$K_n(z_1,z_2)=u_n^{-1}W_n(z_1,z_2)$, where $W_n$ is the kernel of
the $U$-statistics mentioned in Proposition \ref{GU}. Applying
Proposition \ref{GU}, we get
$$
U_n\,\stackrel{P_{n,0}}{\to}\,\zeta\sim \CN(0,1) .
$$
This yields,
\begin{equation}\label{test.1}
E_{n,0}(\psi_{n}^{H})=P_{n,0}(U_n\le -H)=1-\Phi(H)+o(1).
\end{equation}

\subsubsection{Minimax type II error}

By \nref{test.1} we have to verify that
\begin{equation}\label{test.2}
\sup_{f\in\CF(r_n)}E_{n,f}(1-\psi_{n}^{H})=\sup_{f\in\CF(r_n)}P_{n,f}(U_n>H)=\Phi(H-u_n)+o(1).
\end{equation}
For $f=\sum_{l\in\CL}\t_l\phi_l$, let
$$
v_l = \sqrt{n}\t_l,\quad h_n(f)= \frac 12\sum_{l\in\CN}w_{n,l}v_l^2.
$$
\begin{lemma}\label{LT2} Uniformly over $f\in\CF$,
\begin{eqnarray}\label{est}
E_{n,f}U_n &\sim& h_n(f),\label{FF1}\\
\Var_{n,f}U_n&=&1+O(\|f\|^2+\|f\|_4^4). \label{FF2}
\end{eqnarray}
Moreover, uniformly over $f\in\CF$ such that
\begin{equation}
\label{eq:fff} \|f\|=o(1),\quad \|f\|_4=o(1) \quad \text{and} \quad
h_n(f)=O(1),
\end{equation}
the statistics $U_n-h_n(f)$ are asymptotically $\CN(0,1)$, under
$P_{n,f}$-probability.
\end{lemma}

\begin{remark}\label{R3*}
{\rm Using H\"older's inequality and ({\bf A3}) with $p=4+2\delta, \
\delta>0$, we get
$$
\|f\|_4^4\le \|f\|^{a} \|f\|_p^{b},\quad a=2/(1+1/\delta),\quad
b=p/(1+\delta);\quad \|f\|\le\|f\|_p.
$$
Therefore, under ({\bf A3}), Lemma \ref{LT2} yields
\begin{equation}\label{var1}
\sup_{f\in\CF}\Var_{n,f}U_n=O(1)\quad\text{and}\quad\Var_{n,f}U_n=1+O(\|f\|^2+\|f\|^{a})
\end{equation}
uniformly over $f\in\CF$, and $$U_n=h_n(f)+\zeta_n,\
\zeta_n\to\zeta\sim \CN(0,1),$$ uniformly over $f\in\CF$ such that
$h_n(f)=O(1)$ and $\|f\|=o(1)$.}
\end{remark}
{\bf Proof of Lemma \ref{LT2}}. Let the function
$f=n^{-1/2}\sum_{l\in\CL}v_l\phi_l$. Denote $z=(x,t)$ with
$x=f(t)+\xi$, $\xi$ and $t$ are independent, $\xi\sim\CN(0,1)$ and
$t$ is uniformly distributed on $\Delta$. Since the items of the sum
in $U$-statistics are identically distributed and uncorrelated, we
have
$$
E_{n,f}U_n=\frac{n-1}{2}E_{n,f}K_n(z_1,z_2),
$$
where $z_1$ and $z_2$ are independent and distributed as $z$,
\begin{eqnarray*}
E_{n,f}K_n(z_1,z_2)&=&E_{n,f}x_1x_2G_n(t_1,t_2)
=E_n^Tf(t_1)f(t_2)G_n(t_1,t_2)\\
&=&\sum_{l \in \CN}w_{n,l}E_n^T\l(f(t)\phi_l(t)\r)^2=n^{-1}\sum_{l
\in \CN}w_{n,l}v_l^2.
\end{eqnarray*}
Hence, (\ref{FF1}) follows.

Let us now evaluate the variance. Rewrite the $U$-statistics in the
form
\begin{equation}\label{prez}
U_n=  U_{n,0}+U_{n,1}+U_{n,2},
\end{equation}
where $$U_{n,k} = \frac 1n\sum_{1\le i<j\le n}K_{n,k}(z_i,z_j)$$ are
$U$-statistics with the kernels $K_{n,k}(z_1,z_2)$ of the form
\begin{eqnarray*}
K_{n,0}&=&\xi_1\xi_2G_n(t_1,t_2),\quad
K_{n,1}=(\xi_1f(t_2)+\xi_2f(t_1))G_n(t_1,t_2),\\
K_{n,2}&=&f(t_1)f(t_2)G_n(t_1,t_2),\quad G_n(t_1,t_2) = \sum_{l \in
\CN}w_{n,l}\phi_l(t_1)\phi_l(t_2),
\end{eqnarray*}
and the items $U_{n,0}$, $U_{n,1}$ and $U_{n,2}$ are uncorrelated.
Obviously,
\begin{eqnarray*}
E_{n,f}U_{n,0}&=&E_{n,f}U_{n,1}=0,\\
E_{n,f}U_{n,2}&=&\frac{n-1}{2}\sum_{l \in \CN}w_{n,l}\l(\int_\Delta
f(t)\phi_l(t)dt\r)^2\sim h_n(f).
\end{eqnarray*}
Similarly to Proposition \ref{GU},
$$
\Var_{n,f}U_{n,0}\sim \frac 12\int_\Delta\int_\Delta
G_n^2(t_1,t_2)dt_1dt_2=\frac 12\sum_{l \in \CN}w_{n,l}^2=1.
$$
Analogously, by ({\bf A2}) and \nref{cont}, and since
$\max_lw_{n,l}^2=O(1/N)$,
\begin{eqnarray*}
\Var_{n,f}U_{n,1}&\sim &2\int_\Delta\int_\Delta
f^2(t_1)G_n^2(t_1,t_2)dt_1dt_2\\
&=&2 \int_\Delta
\l(f^2(t)\sum_{l\in\CN}w_{n,l}^2\phi_l^2(t)\r)dt=O(\|f\|^2).
\end{eqnarray*}
Next,
$$
\Var_{n,f}U_{n,2}\le \int_\Delta\int_\Delta
f^2(t_1)f^2(t_2)G_n^2(t_1,t_2)dt_1dt_2= A_n.
$$
Let ${\bf{G}}_n$ be the integral operator in $L_2(\Delta)$
associated with the symmetric positively semi-defined kernel
$G_n(t_1,t_2),\ t_1,t_2\in\Delta$, and
$$\|{\bf{G}}_n\|_\infty = \sup_{\|f\|\le
1}\|{\bf{G}}_n f\|=\max_{l\in\CN}w_{n,l}=O(N^{-1/2}).
$$
Observe that, by ({\bf A2}) and \nref{cont},
$$
{\bf{G}}_{n}^{*} =
\sup_{t\in\Delta}\sum_{l\in\CN}w_{n,l}\phi_l^2(t)\le
N\|{\bf{G}}_n\|_\infty,\quad
{\bf{G}}_{n}^{*}\|{\bf{G}}_n\|_\infty=O(1).
$$
We have
\begin{eqnarray*}
A_n&=& \sum_{l\in\CN}w_{n,l}\int_\Delta\int_\Delta
\phi_l(t_1)\phi_l(t_2)f^2(t_1)f^2(t_2)G_n(t_1,t_2)dt_1 dt_2\\
&=& \sum_{l\in\CN}w_{n,l}\langle
f^2\phi_l,{\bf{G}}_{n}(f^2\phi_l)\rangle\le
\|{\bf{G}}_{n}\|_\infty\sum_{l\in\CN}w_{n,l}\|f^2\phi_l\|^2\\
&=&\|{\bf{G}}_{n}\|_\infty\int_\Delta\sum_{l\in\CN}w_{n,l}\phi_l^2(t)f^4(t)dt\\
&\le& \|{\bf{G}}_{n}\|_\infty\sup_{t \in \Delta}
\l(\sum_{l\in\CN}w_{n,l}\phi_l^2(t)\r)\int_\Delta
f^4(t)dt\\&=&\|{\bf{G}}_{n}\|_\infty
{\bf{G}}_{n}^{*}\|f\|^4_4=O(\|f\|^4_4).
\end{eqnarray*}
Hence, (\ref{FF2}) follows.

Using \nref{prez}, and an evaluation similar to the above under
\nref{eq:fff}, we have
$$
U_n-h_n(f)=U_{n,0}+U_{n,1}+U_{n,2}-h_n(f),
$$
where $ U_{n,1}\to 0,\ U_{n,2}-h_n(f)\to 0$, in
$P_{n,f}$-probability. By Proposition \ref{GU}, the statistics
$U_{n,0}$ are asymptotically Gaussian $\CN(0,1)$. This completes the
proof of Lemma \ref{LT2}.
\endproof

\medskip

Let $h_n(f)=O(1)$. Let us now evaluate $\|f\|^2$, $f\in\CF$. We have
$$
\|f\|^2=\sum_{l\in\CL}\t_l^2:=A'_n+B'_n,\quad A'_n =
\sum_{c_l<C/2}\t_l^2, \quad B'_n = \sum_{c_l\ge C/2}\t_l^2.
$$
The second sum is controlled by
$$
B'_n\le 4C^{-2}\sum_{l\in\CL}c_l^2\t_l^2\le 4C^{-2}=o(1).
$$
The first sum is controlled by
\begin{eqnarray*}
A'_n&\le&(4/3)\sum_{l \in \CN}(1-(c_l/C)^2)\t_l^2=(4/3)(w_n/n)\sum_{l \in \CN}w_{n,l}v_n^2\\
&= &(4/3)(w_n/n)h_n(f)=o(h_n(f)),
\end{eqnarray*}
since,  by \nref{cont} and ({\bf B1}), we have $w_n/n=O(
N^{1/2}/n)=o(1)$. Therefore, by \nref{var1}, we have in
$P_{n,f}$-probability,
\begin{equation*}\label{Norm1}
U_n=h_n(f)+\zeta_n,\quad \zeta_n\to\zeta\sim\CN(0,1),
\end{equation*}
uniformly as $h_n(f)=O(1)$.

\begin{lemma}\label{L3}
$$ \inf_{f\in\CF(r_n)}h_n(f)=u_n. $$
\end{lemma}
{\bf Proof of Lemma \ref{L3}} It follows using general convexity
arguments (see \cite{I.93}, Lemma 11 of \cite{IK.07}, Proposition
4.1 of \cite{IS.02}).
\endproof

\medskip

Let us now evaluate type II errors for a sequence
$f=f_n\in\CF(r_n)$. First, let $h_n(f_n)\to\infty$. Applying Lemmas
\ref{LT2}, \ref{L3} and \nref{var1}, we  have
\begin{eqnarray*}
E_{n,f}(1-\psi_{n}^{H})&=&P_{n,f}(U_n\le H)=P_{n,f}(E_{n,f}-U_n\ge
E_{n,f}-H)\\
&\le&\Var_{n,f}(U_n)/(E_{n,f}-H)^2=o(1).
\end{eqnarray*}
Let $h_n(f_n)=O(1)$ (by Lemma \ref{L3} this is only possible for
$u_n=O(1)$). Applying Lemmas \ref{LT2}, \ref{L3} and \nref{var1}
once again, we have
\begin{eqnarray*}
E_{n,f}(1-\psi_{n}^{H})&=&P_{n,f}(U_n\le H)=P_{n,f}(E_{n,f}-U_n\ge
E_{n,f}-H)\\
&=&P_{n,f}(\zeta_n\ge h_n(f)-H+o(1))=\Phi(H-h_n(f))+o(1).
\end{eqnarray*}
Therefore,
$$
\sup_{f\in\CF(r_n)}E_{n,f}(1-\psi_{n}^{H})=\Phi(H-\inf_{f\in\CF(r_n)}h_n(f))+o(1)=
\Phi(H-u_n)+o(1).
$$
This yields \nref{test.2}. Hence, Theorem \ref{T1} (2) follows.
\endproof

\medskip

This completes the proof of Theorem \ref{T1}.

\subsection{Upper bounds for Theorem \protect{\ref{T0}}}
Observe that the kernel \nref{test1} is of the form \nref{test} with
coefficients $$w_{l,n}=w_n=\sqrt{2/N},\quad l\in\CN.$$ Hence,
Proposition \ref{GU} is applicable to the $U$-statistics $U_n$ with
kernel \nref{test1} and yields asymptotic normality $\CN(0,1)$ of
$U_n$ under $P_{n,0}$. Thus, we get \nref{test.1}. Analogously, we
obtain Lemma \ref{LT2} with
$$
h_n(f)=\frac{n}{\sqrt{2N}}\sum_{l \in \CN}\t_l^2.
$$
If $h_n(f)=O(1), \
f\in\CF$, then $\|f\|=o(1)$. In fact,
$$
\|f\|^2=\sum_{l\in\CL}\t_l^2\le \sum_{l \in
\CN}\t_l^2+C^{-2}\sum_{c_l\ge C}c_l^2\t_l^2\le
\frac{\sqrt{2N}}{n}h_n(f)+C^{-2}=o(1).
$$
These yield \nref{test.2} for $f\in\CF$ such that $h_n(f)=O(1)$. If
$h_n(f)\to\infty$, then it follows from Chebyshev's inequality and
the boundness of the variances that $P_{n,f}(U_n\ge H)\to 0$ for
$H<ch_n(f),\ c\in (0,1)$. Hence, Theorem \ref{T0} (2) follows.
\endproof

\medskip
This completes the proof of Theorem \ref{T0}.

\end{document}